\documentclass{gtart_h}  

\def\ifplaintex{\expandafter\ifx\csname documentclass\endcsname\relax}

\ifplaintex 
\else
\fi

\expandafter\ifx\csname beginpicture\endcsname\relax
\expandafter\ifx\csname documentclass\endcsname\relax
\input pictex \else
\input prepictex \input pictex \input postpictex \fi\fi

\def\gt{{\mathsurround=0pt\it $\cal G\mskip-2mu$eometry \&\ 
$\cal T\!\!$opology}}        

\def\gtp{{\mathsurround=0pt\it $\cal G\mskip-2mu$eometry \&\ 
$\cal T\!\!$opology $\cal P\!$ublications}}  

\def\lognumber#1{\def\thelognumber{#1}}
\def\volumenumber#1{\def\thevolumenumber{#1}}
\def\papernumber#1{\def\thepapernumber{#1}}
\def\volumeyear#1{\def\thevolumeyear{#1}}

\def\pagenumbers#1#2{\def\startpage{#1}\def\finishpage{#2}}
\def\published#1{\def\publishdate{#1}}
\def\proposed#1{\def\theproposer{#1}}
\def\seconded#1{\def\theseconders{#1}}
\def\received#1{\def\receiveddate{#1}}
\def\revised#1{\def\reviseddate{#1}}
\def\accepted#1{\def\accepteddate{#1}}
\def\asciititle#1{\def\theasciititle{#1}}

\def\asciiaddress#1{\def\theasciiaddress{#1}}
\def\asciiemail#1{\def\theasciiemail{#1}}

\long\def\asciiabstract#1{\long\def\theasciiabstract{#1}}

\def\shortauthors#1{\def\theshortauthors{#1}}

\let\\\par\let\thelognumber\relax
\let\thevolumenumber\relax\let\thepapernumber\relax
\let\thevolumeyear\relax\let\thesamplenumber\relax\let\startpage\relax
\let\finishpage\relax\let\publishdate\relax\let\receiveddate\relax
\let\reviseddate\relax\let\accepteddate\relax\let\theasciititle\relax
\let\theasciiauthors\relax\let\theasciiaddress\relax
\let\theasciiabstract\relax
\let\theasciiemail\relax\let\theshortauthors\relax\let\theshorttitle\relax

\long\def\maketitlep{   

\count0=\startpage

\gt\hfill      
\beginpicture
\setcoordinatesystem units <0.33truein, 0.33truein> point at 2.2 0.9
\setplotsymbol ({$\cal G$})
\plotsymbolspacing=9truept
\circulararc 315 degrees from 0 1 center at 0 0
\setplotsymbol ({$\cal T$})
\circulararc 315 degrees from 1 -1 center at 1 0
\endpicture
\break
{\small\ifx\thesamplenumber\relax 
Volume \else Sample
\fi\thevolumenumber\ (\thevolumeyear)
\startpage--\finishpage\nl
Published: \publishdate}
\vglue 0.5truein plus 0.4fil minus 0.1truein

{\parskip=0pt\leftskip 0pt plus 1fil\def\\{\par\smallskip}{\ifplaintex\large
\else\Large\fi\bf\thetitle}\par\medskip}   

\vglue 0pt plus 0.1fil 

{\parskip=0pt\leftskip 0pt plus 1fil\def\\{\par}{\sc\theauthors}
\par\medskip}

\vglue 0pt plus 0.1fil 

{\small\parskip=0pt\let\newline\\
{\leftskip 0pt plus 1fil\def\\{\par}{\sl\theaddress}\par}
\expandafter\ifx\theemail\relax    
\relax\else\vglue 5pt plus 0.02fil minus 2pt\def\\{\stdspace{\rm 
and}\stdspace} 
\cl{Email:\stdspace\tt\theemail}\fi
\ifx\theurl\relax                  
\relax\else\vglue 5pt plus 0.02fil minus 2pt\def\\{\stdspace{\rm 
and}\stdspace}
\cl{URL:\stdspace\tt\theurl}\fi\par}

\vglue 7pt plus 0.3fil minus 3pt

{\bf Abstract}
\vglue 5pt plus 0.1fil minus 2pt

\theabstract

\vglue 7pt plus 0.3fil minus 3pt

{\bf AMS Classification numbers}\quad Primary:\quad \theprimaryclass

Secondary:\quad \thesecondaryclass

\vglue 5pt plus 0.3fil minus 2pt

{\bf Keywords:}\quad \thekeywords

\vglue 10pt plus 0.5fil minus 5pt

{\small  Proposed: \theproposer\hfill Received: \receiveddate\nl
Seconded: \theseconders\hfill 
\ifx\reviseddate\relax                         
Accepted: \accepteddate                        
\else
Revised: \reviseddate                          
\fi}
\eject
}       

\let\maketitlepage\maketitlep
\let\maketitle\maketitlepage

\font\phead=cmsl9 scaled 950
\font=cmsl9 scaled 1050
\font\pnum=cmbx10 scaled 913
\font\lnum=cmbx10 
\font\pfoot=cmsl9 scaled 950
\font=cmsl9 scaled 1050
\ifplaintex
\headline{\vbox to 0pt{\vskip -4.5mm\line{\small\phead\ifnum
\count0=\startpage ISSN 1364-0380 (on line)
1465-3060 (printed) \hfill {\pnum\folio}\else\ifodd\count0\def\\{ }%
\ifx\theshorttitle\relax\thetitle\else\theshorttitle\fi\hfill{\pnum\folio}
\else\def\\{ and }{\pnum\folio}\hfill\ifx\theshortauthors\relax\theauthors
\else\theshortauthors\fi\fi\fi}\vss}}
\footline{\vbox to 0pt{\vglue 0mm\line{\small\pfoot\ifnum\count0=\startpage
\copyright\ \gtp\hfill\else
\gt, Volume \thevolumenumber\ (\thevolumeyear)\hfill\fi}\vss
}}
\else
\makeatletter
\def\@oddhead{{\small\lhead\ifnum\count0=\startpage ISSN 1364-0380 (on line)
1465-3060 (printed) \hfill {\lnum\number\count0}\else\ifodd\count0
\def\\{ }\ifx\theshorttitle\relax \thetitle \else\theshorttitle\fi\hfill
{\lnum\number\count0}\else\def\\{ and }{\lnum\number\count0}
\hfill\ifx\theshortauthors\relax 
\theauthors\else\theshortauthors\fi\fi\fi}}\def\@evenhead{\@oddhead}
\def\@oddfoot{\small\lfoot\ifnum\count0=\startpage\copyright\ \gtp\hfill\else
\gt, Volume \thevolumenumber\ (\thevolumeyear)\hfill\fi}
\def\@evenfoot{\@oddfoot}
\makeatother
\fi

\newwrite\gtoutfile
\long\gdef\makeheadfile{  
{\def\\{, }\def\s{ }
\immediate\openout\gtoutfile head.xxx
\immediate\write\gtoutfile{Proxy-for: \ifx\theasciiauthors\relax
\theauthors\else\theasciiauthors\fi\s<\ifx\theasciiemail\relax\theemail\else\theasciiemail\fi>}
\immediate\write\gtoutfile{\noexpand\\}
\immediate\write\gtoutfile{Authors: \ifx\theasciiauthors\relax
\theauthors\else\theasciiauthors\fi}
{\def\\{ }\immediate\write\gtoutfile{Title: \ifx\theasciititle\relax
\thetitle\else\theasciititle\fi}}
\immediate\write\gtoutfile{Subj-class: GT or SG or MG etc}
\immediate\write\gtoutfile{MSC-class: \theprimaryclass\ifx\thesecondaryclass\relax\else, \thesecondaryclass\fi}
\immediate\write\gtoutfile{Journal-ref: Geom. Topol. \thevolumenumber
(\thevolumeyear) \startpage-\finishpage}
\immediate\write\gtoutfile{Comments: Published by Geometry and Topology at}
\immediate\write\gtoutfile{\s\s http://www.maths.warwick.ac.uk/gt/GTVol\thevolumenumber/paper\thepapernumber.abs.html}
\immediate\write\gtoutfile{\noexpand\\}
\immediate\write\gtoutfile{}
\ifx\theasciiabstract\relax
\immediate\write\gtoutfile{\theabstract}\else
\immediate\write\gtoutfile{\theasciiabstract}\fi
\immediate\write\gtoutfile{}
\immediate\write\gtoutfile{\noexpand\\}
\immediate\write\gtoutfile{}
\immediate\closeout\gtoutfile}}  

\def\maketitlepage{\maketitlep\makeheadfile}
\let\maketitle\maketitlepage

\lognumber{591}

\received{5 April 2005}
\volumenumber{9}\papernumber{34}\volumeyear{2005}
\pagenumbers{1501}{1538}   
\revised{25 July 2005}
\published{8 August 2005}
\accepted{24 June 2005}
\proposed{Walter Neumann}
\seconded{Martin Bridson, Benson Farb}

\usepackage{amsmath,amssymb}

\numberwithin{equation}{subsection}
\theoremstyle{plain}
\newtheorem{theorem}[equation]{Theorem}
\newtheorem{thm}[equation]{Theorem}

\newtheorem{prop}[equation]{Proposition}
\newtheorem{lem}[equation]{Lemma}
\newtheorem{lemma}[equation]{Lemma}

\newtheorem{corollary}[equation]{Corollary}
\newtheorem{cor}[equation]{Corollary}

\theoremstyle{remark}

\theoremstyle{definition}

\newtheorem{defn}[equation]{Definition}

\newtheorem{example}[equation]{Example}

\newcommand{\nbd}[2]{\mathcal{N}_{#2} ({#1}) } 
\newcommand{\bignbd}[2]{\mathcal{N}_{#2} \bigl({#1}\bigr) }

\newcommand{\ra}{\rightarrow}

\newcommand{\restr}{\mbox{\Large \(|\)\normalsize}}

\newcommand{\A}{{\mathcal A}}

\newcommand{\D}{{\mathcal D}}
\newcommand{\E}{\mathbb E}
\newcommand{\F}{{\mathcal F}}

\newcommand{\N}{\mathbb N}
\renewcommand{\P}{\mathbb P}

\newcommand{\Z}{\mathbb Z}

\newcommand{\Hyp}{\mathbb H}
\newcommand{\Set}[1]{\mathcal{#1}}

\DeclareMathOperator{\diam}{diam}
\renewcommand{\flat}{\operatorname{Flat}}

\newcommand{\isom}{\operatorname{Isom}}

\newcommand{\acts}{\curvearrowright}
\DeclareMathOperator{\CAT}{CAT}
\DeclareMathOperator{\Cone}{Cone}
\DeclareMathOperator{\Isom}{Isom}
\DeclareMathOperator{\Stab}{Stab}

\def\D{\partial}
\newcommand{\boundary}{\partial}
\newcommand{\al}{\alpha}
\def\de{\delta}

\def\eps{\epsilon}
\def\ga{\gamma}
\def\Ga{\Gamma}

\newcommand{\oa}{\overrightarrow}

\def\om{\omega}

\def\ra{\rightarrow}

\def\Si{\Sigma}

\def\th{\theta}
\def\ulim{\mathop{\hbox{\textup{$\om$-lim}}}}
\def\geo{\partial_{\infty}}

\newcommand{\defeq}{\mathbin{:=}}
\newcommand{\rdefeq}{\mathbin{=:}}

\newcommand{\inclusion}{\hookrightarrow}

\newcommand{\ol}{\overline}

\def\cangle{\widetilde\angle} 
\def\tangle{\angle_{T}}       
\def\Td{d_T}                  
\def\tits{\partial_T}         

\newcommand{\set}[2]{\{\, {#1} \mid {#2} \,\}}
\newcommand{\bigset}[2]{ \bigl\{ \, {#1} \bigm| {#2} \, \bigr\} }

\def\RomanianComma#1{\setbox0=\hbox{#1}{\ooalign{\hidewidth
    \lower1.2ex\hbox{$\mspace{1mu}^{,}$}\hidewidth\crcr\unhbox0}}}
\newcommand{\Drutu}{Dru{\RomanianComma{t}u}}

\begin{document}

\title{Hadamard spaces with isolated flats}
\asciititle{Hadamard spaces with isolated flats, 
with an appendix written jointly with Mohamad Hindawi}

\authors{G Christopher Hruska\\Bruce Kleiner}
\shortauthors{G\,C Hruska and B Kleiner}

\address{{\rm With an appendix written jointly with {\sc Mohamad Hindawi}}\\
\\\smallskip\\{\rm GCH:}\qua
Department of Mathematics, 
University of Chicago\\5734 S~University 
Ave, Chicago, IL 60637-1514, USA\\\smallskip\\{\rm BK:}\qua Department of Mathematics,
University of Michigan\\Ann Arbor, MI 48109-1109, 
USA\\\smallskip\\{\rm MH:}\qua Department of Mathematics, University of 
Pennsylvania\\Philadelphia, PA 19104-6395, USA\\\smallskip\\{\rm Email:}\qua 
\tt\mailto{chruska@math.uchicago.edu}, 
\mailto{bkleiner@umich.edu}\\\mailto{mhindawi@math.upenn.edu}}

\asciiaddress{GCH:
Department of Mathematics, 
University of Chicago\\5734 S University 
Ave, Chicago, IL 60637-1514, USA\\BK: Department of Mathematics,
University of Michigan\\Ann Arbor, MI 48109-1109, 
USA\\MH: Department of Mathematics, University of 
Pennsylvania\\Philadelphia, PA 19104-6395, USA}

\asciiemail{chruska@math.uchicago.edu, bkleiner@umich.edu, 
mhindawi@math.upenn.edu}

\begin{abstract}
We explore the geometry of nonpositively curved spaces with isolated flats,
and its consequences for groups that act properly discontinuously,
cocompactly, and isometrically on such spaces.
We prove that the geometric boundary of the space is an invariant
of the group up to equivariant homeomorphism.
We also prove that any such group is relatively hyperbolic,
biautomatic, and satisfies the Tits Alternative.
The main step in establishing these results
is a characterization of spaces with isolated flats as
relatively hyperbolic with respect to flats.
Finally we show that a $\CAT(0)$ space has isolated flats
if and only if its Tits boundary is a disjoint union
of isolated points and standard Euclidean spheres.

In an appendix written jointly with Hindawi, we extend many of the
results of this article to a more general setting in which the
isolated subspaces are not required to be flats.
\end{abstract}

\asciiabstract{%
We explore the geometry of nonpositively curved spaces with isolated
flats, and its consequences for groups that act properly
discontinuously, cocompactly, and isometrically on such spaces.  We
prove that the geometric boundary of the space is an invariant of the
group up to equivariant homeomorphism.  We also prove that any such
group is relatively hyperbolic, biautomatic, and satisfies the Tits
Alternative.  The main step in establishing these results is a
characterization of spaces with isolated flats as relatively
hyperbolic with respect to flats.  Finally we show that a CAT(0) space
has isolated flats if and only if its Tits boundary is a disjoint
union of isolated points and standard Euclidean spheres.

In an appendix written jointly with Hindawi, we extend many of the
results of this article to a more general setting in which the
isolated subspaces are not required to be flats.
}



\maketitlepage

\section{Introduction}
\label{sec:Introduction}

In this article, we explore the large scale geometry of $\CAT(0)$
spaces with isolated flats and its implications for groups that act
geometrically, ie, properly discontinuously, cocompactly,
and isometrically, on such spaces.
Spaces with isolated flats have many features in common with
Gromov-hyperbolic spaces and can be viewed as the nonhyperbolic
$\CAT(0)$ spaces that are closest to being hyperbolic.

Throughout this article, a \emph{$k$--flat} is an isometrically embedded
copy of Euclidean space $\E^k$ for $k\ge 2$.  In particular
a geodesic line is not considered to be a flat.
Let $\flat(X)$ denote the space of all flats in~$X$ with the topology of
Hausdorff convergence on bounded sets
(see Definition~\ref{def:PointedHausdorff} for details).
A $\CAT(0)$ space with a geometric group action
has \emph{isolated flats} if it contains an equivariant collection~$\F$
of flats such that $\F$ is closed and isolated in $\flat(X)$ and
each flat $F \subseteq X$ is contained
in a uniformly bounded tubular neighborhood of some $F' \in \F$.
As with the notion of Gromov-hyperbolicity, the notion of isolated flats
can be characterized in many equivalent ways.

\subsection{Examples}
\label{subsec:Examples}

The prototypical example of a $\CAT(0)$ space with isolated flats is the
truncated hyperbolic space associated to a finite volume cusped
hyperbolic manifold~$M$.  Such a space is obtained from hyperbolic
space~$\Hyp^n$ by removing an equivariant collection of disjoint
open horoballs that are the lifts of the cusps of~$M$
and endowing the resulting space with the induced length metric.
The truncated space is nonpositively curved,%
\footnote{
   This fact is a special case of a result stated by Gromov in
   \cite[\S 2.2]{Gromov81}.
   A proof of Gromov's theorem was provided by Alexander--Berg--Bishop
   in \cite{AlexanderBergBishop93}.
   For a proof tailored to this special case, see \cite[II.11.27]{BH99}.
}
and its only flats are the boundaries of the deleted horoballs.
It is easy to see that the truncated space has isolated flats.
Furthermore, the fundamental group of~$M$ acts cocompactly on the
truncated space.

Other examples of groups that act geometrically on $\CAT(0)$ spaces with
isolated flats include the following:
\begin{itemize}
  \item Geometrically finite Kleinian groups \cite{Bowditch93}.
  \item Fundamental groups of compact manifolds obtained by gluing
  finite volume hyperbolic manifolds along cusps
  (Heintze, see \cite[\S 3.1]{Gromov81}).
  \item Fundamental groups of compact manifolds obtained by the cusp closing
  construction of Thurston--Schroeder \cite{Schroeder89}.
  \item Limit groups (also known as $\omega$--residually free groups),
  which arise in the study of equations over free groups
  \cite{AlibegovicBestvinaLimit}.
\end{itemize}

Examples in the $2$--dimensional setting include the fundamental group 
of any compact nonpositively curved $2$--complex whose $2$--cells are
isometric to regular Euclidean hexagons \cite{BallmannBrin94,Wise96}.
Ballmann--Brin showed that such $2$--complexes exist in abundance
and can be constructed with arbitrary local data
\cite{BallmannBrin94}.
For instance, for each simplicial graph~$L$
there is a $\CAT(0)$ hexagonal $2$--complex~$X$ such that 
the link of every vertex in~$X$ is isomorphic to the graph~$L$
(\cite{Moussong88}, see also \cite{Haglund91} and \cite{Benakli94}).

As indicated in \cite{KapovichLeeb95}, the notion of isolated flats
has a natural generalization where the family~$\F$ is a collection of
closed convex subspaces rather than a collection of flats.
In an appendix by Hindawi, Hruska, and Kleiner, we extend the results
of this article to this more general setting.

In particular, this generalization includes the universal covers of
compact $3$--manifolds whose geometric decomposition contains
at least one hyperbolic component (see \cite{KapovichLeeb95}),
and the universal covers of closed, real analytic, nonpositively curved
$4$--manifolds whose Tits boundary does not contain a nonstandard component
(see \cite{Hindawi05} for details).

\subsection{Main results}
\label{subsec:MainResults}

The following theorem proves that several other conditions are equivalent
to having isolated flats.  Further equivalent geometric notions
are discussed in Theorem~\ref{thm:GeometricEquivalence} below.

\begin{thm}\label{thm:MainEquivalence}
Let $X$ be a $\CAT(0)$ space and $\Gamma$ a group acting geometrically
on~$X$.  The following are equivalent.
   \begin{enumerate}
   \item \label{item:MainIsolated}
   $X$ has isolated flats.
   \item \label{item:MainTits}
   Each component of the Tits boundary $\tits X$ is either an
   isolated point or a standard Euclidean sphere.
   \item \label{item:MainCone}
   $X$ is a relatively hyperbolic space with respect to a family of
   flats~$\F$.
   \item \label{item:MainRelHyp}
   $\Gamma$ is a relatively hyperbolic group with respect to a collection of
   virtually abelian subgroups of rank at least two.
   \end{enumerate}
\end{thm}

Several different nomenclatures exist in the literature for
relative hyperbolicity.
For groups, we use the terminology of Bowditch throughout.
In Farb's terminology, the groups we call ``relatively hyperbolic''
are called ``relatively hyperbolic with Bounded Coset Penetration.''
For metric spaces, the property we call ``relatively hyperbolic''
was introduced by \Drutu--Sapir, although they used the
term ``asymptotically tree-graded'' for such spaces.
\Drutu--Sapir proved the equivalence of the metric and group theoretic
notions of relative hyperbolicity for a finitely generated group
with the word metric \cite{DrutuSapir05}.

The implications
(\ref{item:MainCone})~$\Longleftrightarrow$~(\ref{item:MainRelHyp})
and (\ref{item:MainCone})~$\Longrightarrow$~(\ref{item:MainIsolated})
follow in a straightforward fashion from work of
\Drutu--Osin--Sapir \cite{DrutuSapir05}.
A large part of this article consists of establishing the remaining 
implications
(\ref{item:MainIsolated})~$\Longrightarrow$~(\ref{item:MainCone})
and (\ref{item:MainIsolated})~$\Longleftrightarrow$~(\ref{item:MainTits}).

A result analogous to
(\ref{item:MainIsolated})~$\Longrightarrow$~(\ref{item:MainCone})
is established by Kapovich--Leeb in \cite{KapovichLeeb95}
using a more restrictive notion of isolated flats, discussed in more detail
in the next subsection.  This notion, while sufficiently general for
their purposes, is tailored to make their proof go smoothly.
A major goal of this article is to generalize Kapovich--Leeb's results
to a more natural class of spaces.
A portion of our proof can be viewed as a
streamlined version of Kapovich--Leeb's proof.

Theorem~\ref{thm:MainEquivalence}
has the following consequences for groups acting on spaces with
isolated flats, using some existing results from the literature.

\begin{thm}
\label{thm:Applications}
Let $\Gamma$ act geometrically on a
$\CAT(0)$ space~$X$ with isolated flats.  Then $X$ and $\Gamma$ have the
following properties.
   \begin{enumerate}
   \item \label{item:Quasiflat}
   Quasi-isometries of~$X$ map maximal flats to maximal flats.
   \item \label{item:UndistortedQC}
   A finitely generated
   subgroup $H \le \Gamma$ is undistorted if and only if it is
   quasiconvex \textup{(}with respect to the $\CAT(0)$ action\textup{)}.
   \item \label{item:BoundaryWellDefined}
   The geometric boundary $\boundary X$ is a group invariant of~$\Gamma$.
   \item \label{item:TitsAlt}
   $\Gamma$ satisfies the Strong Tits Alternative.  In other words,
   every subgroup of\/~$\Gamma$ either is virtually abelian or contains
   a free subgroup of rank two.
   \item \label{item:Biautomatic}
   $\Gamma$ is biautomatic.
   \end{enumerate}
\end{thm}

Property~(\ref{item:Quasiflat}) is a direct consequence of
Theorem~\ref{thm:MainEquivalence} and a result of \Drutu--Sapir
\cite[Proposition~5.4]{DrutuSapir05}.
Properties (\ref{item:UndistortedQC}) and~(\ref{item:BoundaryWellDefined})
are consequences of Theorem~\ref{thm:MainEquivalence} and
results proved by Hruska in
\cite{HruskaGeometric}.
Property~(\ref{item:TitsAlt}) follows from
Theorem~\ref{thm:MainEquivalence} together with work of Gromov and Bowditch
\cite{Gromov87,BowditchRelHyp}.
Property~(\ref{item:Biautomatic}) is proved using
Theorem~\ref{thm:MainEquivalence} and a result of Rebbechi
\cite{Rebbechi01}.

In addition,
Chatterji--Ruane have used Theorem~\ref{thm:MainEquivalence} 
to show that a group acting geometrically on a space with isolated flats
satisfies the property of Rapid Decay as well as the Baum--Connes
conjecture \cite{ChatterjiRuane05} (cf \cite{DrutuSapirRD05}).
Groves has also used Theorem~\ref{thm:MainEquivalence}
in his proof that toral, torsion free groups acting geometrically
on spaces with isolated flats are Hopfian \cite{GrovesHopfian}.

In proving Theorem~\ref{thm:MainEquivalence},
we found it useful to clarify the notion of isolated flats
by exploring various geometric formulations and proving their equivalence.
Indeed, we consider this exploration to be the most novel contribution
of this article.  Each of these formulations has advantages in certain situations.  We explain these various formulations below.

Let $X$ be a $\CAT(0)$ space and $\Gamma$ a group acting geometrically
on~$X$.  We say that $X$ has (IF1), or 
\emph{isolated flats in the first sense},
if it satisfies the definition of
isolated flats given above.
A flat in~$X$ is \emph{maximal} if it is not contained in a finite
tubular neighborhood of any higher dimensional flat.
The space~$X$ has \emph{thin parallel sets}
if for each geodesic line $\gamma$ in a maximal flat~$F$
the parallel set $\P(\gamma)$ lies in a finite tubular neighborhood of~$F$.
The space has \emph{uniformly thin parallel sets}
if there is a uniform bound on the thickness of these tubular neighborhoods.
The space~$X$ has \emph{slim parallel sets}
if for each geodesic $\gamma$ in a maximal flat~$F$, the Tits boundary
of the parallel set $\P(\gamma)$ is equal to the Tits boundary of~$F$.

We say that $X$ has (IF2), or \emph{isolated flats in the second sense},
if there is a $\Gamma$--invariant set~$\F$ of flats in~$X$ such that
the following properties hold.
\begin{enumerate}
  \item \label{item:AllFlats}
  There is a constant $D < \infty$ so that each flat in~$X$
  lies in a $D$--tubular neighborhood of some flat $F \in \F$.
  \item \label{item:ControlledInt}
  For each positive $\rho< \infty$ there is a constant
  $\kappa = \kappa(\rho) < \infty$ so that for any two distinct flats
  $F,F' \in \F$ we have
  \[
     \diam \bigl( \nbd{F}{\rho} \cap \nbd{F'}{\rho} \bigr) < \kappa.
  \]
\end{enumerate}

\begin{thm}\label{thm:GeometricEquivalence}
Let $X$ be a $\CAT(0)$ space and $\Gamma$ a group acting geometrically
on~$X$.  The following are equivalent.
  \begin{enumerate}
  \item \label{item:IF1} $X$ has \textup{(IF1)}.
  \item \label{item:UniformlyThin} $X$ has uniformly thin parallel sets.
  \item \label{item:Thin} $X$ has thin parallel sets.
  \item \label{item:Slim} $X$ has slim parallel sets.
  \item \label{item:IF2} $X$ has \textup{(IF2)}.
  \end{enumerate}
\end{thm}

\noindent
Of these formulations, (IF1) seems the most natural
and also appears to be the most easily verified in practice.
For instance, (IF1) is evidently satisfied for the
cusp gluing manifolds of Heintze and the cusp closing manifolds of
Thurston--Schroeder.
The notion of slim parallel sets plays a crucial role in the
proof of (\ref{item:MainTits})~$\Longrightarrow$~(\ref{item:MainIsolated})
from Theorem~\ref{thm:MainEquivalence}.
The notion of thin parallel sets is important in establishing
the equivalence of the remaining formulations.
The formulation (IF2), while rather technical in appearance,
seems to be the strongest and most widely 
applicable formulation, giving precise control over the coarse intersections
of any pair of maximal flats.
This formulation is used throughout this paper and has 
previously been used in the literature as a definition of isolated
flats (see \cite{Hruska2ComplexIFP,HruskaGeometric}).

\subsection{Historical background}
\label{subsec:History}

$\CAT(0)$ spaces with isolated flats were first considered
by Kapovich--Leeb and Wise, independently.
Kapovich--Leeb \cite{KapovichLeeb95} study a class of $\CAT(0)$ spaces
in which
the maximal flats are disjoint and separated by regions of strict
negative curvature, ie, the separating regions are
locally $\CAT(-1)$.
It is clear that such a space satisfies (IF1).
As mentioned above, they prove that such a space is relatively
hyperbolic, although they use a different terminology,
and they conclude that quasi-isometries of such spaces
coarsely preserve the set of all flats.
The latter result is a key step in their quasi-isometry
classification of nongeometric Haken $3$--manifold groups.

The Flat Closing Problem asks the following:  If $\Gamma$ acts geometrically
on a $\CAT(0)$ space~$X$, is it true that either $\Gamma$ is word
hyperbolic or $\Gamma$ contains a subgroup isomorphic to $\Z \times \Z$?
If $X$ has isolated flats, an easy geometric argument
solves the Flat Closing Problem in the affirmative.
Such an argument was discovered by Ballmann--Brin in the context
of $\CAT(0)$ hexagonal $2$--complexes \cite{BallmannBrin94}
and independently by Wise for arbitrary $\CAT(0)$ spaces satisfying the
weaker formulation (IF1) of isolated flats \cite[Proposition~4.0.4]{Wise96}.

The Flat Plane Theorem states that a proper, cocompact $\CAT(0)$ space
either is $\delta$--hyperbolic or contains a flat plane
\cite{Eberlein73,Gromov87,Heber87,Bridson95}.
In the $2$--dimensional setting, Wise proved an analogous
Flat Triplane Theorem, which states that a proper, cocompact
$\CAT(0)$ $2$--complex either has isolated flats or contains an
isometrically embedded triplane.  A \emph{triplane} is the space formed by
gluing three Euclidean halfplanes isometrically along their boundary lines.
A proof, due to Wise, of the Flat Triplane Theorem first appeared in 
\cite{Hruska2ComplexIFP},
however the ideas are implicit in Wise's article
\cite{WiseFigure8}, which has been circulated since 1998.
Some of the ideas in the proof of Theorem~\ref{thm:GeometricEquivalence}
have combinatorial analogues in the proof of the Flat Triplane Theorem,
although the $2$--dimensional situation is substantially simpler
than the general case.

The question of whether (or when) the geometric boundary depends
only on the group $\Gamma$ was raised in
\cite[\S 6.{\(\text{B}_4\)}]{Gromov93}; in the Gromov hyperbolic
case the invariance follows from the stability of quasi-geodesics.
Croke--Kleiner \cite{CrokeKleiner00,CrokeKleiner02}
found examples which showed that this is
not always the case, and then analyzed, for a class of examples including
nonpositively curved graph manifolds, exactly what geometric
structure determines the geometric boundary up to equivariant
homeomorphism.  We note that Buyalo, using a different geometric idea,
later found examples showing
that the equivariant homeomorphism type of the boundary
is not a group invariant \cite{Buyalo98}.
At the beginning of Croke--Kleiner's project,
Kleiner studied, in unpublished work from 1997, the simpler
case of spaces with isolated flats;
he proved (a reformulation of) the implications 
(\ref{item:MainIsolated})~$\Longrightarrow$ (\ref{item:MainTits}),
(\ref{item:MainCone}) of Theorem~\ref{thm:MainEquivalence},
using the strong hypothesis (IF2),
and showed that equivariant quasi-isometries map
geodesic segments close to geodesic segments and induce
equivariant boundary homeomorphisms.

Hruska studied isolated flats in the $2$--dimensional setting
in \cite{Hruska2ComplexIFP}.  He showed that in that setting,
isolated flats is equivalent to the Relatively Thin Triangle Property,
which states that each geodesic triangle is ``thin relative to a flat,''
and also to the Relative Fellow Traveller Property, which states
that quasigeodesics with common endpoints track close together
``relative to flats.''
The Relative Fellow Traveller Property generalizes a phenomenon
discovered by Epstein in the setting of geometrically finite Kleinian
groups \cite[Theorem 11.3.1]{ECHLPT92}.
In \cite{HruskaGeometric}
Hruska then established parts (\ref{item:Quasiflat}),
(\ref{item:UndistortedQC}), and (\ref{item:BoundaryWellDefined}) 
of Theorem~\ref{thm:Applications}
for $\CAT(0)$ spaces satisfying both (IF2) and the Relative Fellow
Traveller Property.

\subsection{Summary of the sections}
\label{subsec:Summary}

Section~\ref{sec:Background} contains general facts about
the geometry of $\CAT(0)$ spaces and asymptotic cones.
Subsection~\ref{subsec:Notation} is a brief review of facts
from the literature, serving to establish the terminology
and notation we use throughout.
In the remaining subsections, we establish preliminary
results that are not specific to the isolated flats setting.
In Subsection~\ref{subsec:SmallDeficit}, we prove several results 
about triangles in which
one or more vertex angles have a small angular deficit.
In Subsection~\ref{subsec:PeriodicEuclidean},
we prove a lemma characterizing periodic flats whose boundary
sphere is isolated in the Tits boundary of the $\CAT(0)$ space.
In Subsection~\ref{subsec:TriangleUltralimit},
we prove a result about ultralimits of triangles in an ultralimit of 
$\CAT(0)$ spaces.

The goal of
Section~\ref{sec:IsolatedTreeGraded} is to prove the implication
(\ref{item:MainIsolated})~$\Longrightarrow$~(\ref{item:MainCone})
from Theorem~\ref{thm:MainEquivalence}.
The definition of isolated flats used in this section is the strong
formulation (IF2).
Subsection~\ref{subsec:Isolated} is
a review of basic properties of spaces
with isolated flats; in particular, we prove 
that maximal flats are periodic (this result is due to Ballmann--Brin
and Wise, independently).
In Subsection~\ref{subsec:NearAFlat}
we prove Proposition~\ref{prop:TriangleNearFlat},
which states that, in the presence of isolated flats,
a large nondegenerate triangle with small angular deficit
must lie close to a flat (cf Proposition~4.2 of 
\cite{KapovichLeeb95}).
This proposition can be viewed as a generalization
of the well-known result that a triangle with zero angle deficit
bounds a flat Euclidean region.
In Subsection~\ref{subsec:IsolatedIsTreeGraded}
we study the geometry of the asymptotic cones of a space with isolated
flats.  A key result is Proposition~\ref{prop:DirectionComponents},
which determines the space of directions at a point in an asymptotic cone.
The space of directions consists of a disjoint union of spherical components
corresponding to flats and a discrete set of isolated directions. 
Finally we prove Theorem~\ref{thm:TreeGraded}, which establishes
the implication
(\ref{item:MainIsolated})~$\Longrightarrow$~(\ref{item:MainCone})
from Theorem~\ref{thm:MainEquivalence}.  As mentioned above,
the proof of
(\ref{item:MainCone})~$\Longrightarrow$~(\ref{item:MainIsolated})
is immediate from work of \Drutu--Sapir.

In Section~\ref{sec:Applications}, we focus on applications
of Theorem~\ref{thm:TreeGraded}.  We briefly explain the equivalence
(\ref{item:MainCone})~$\Longleftrightarrow$~(\ref{item:MainRelHyp})
of Theorem~\ref{thm:MainEquivalence}
and explain the various parts of Theorem~\ref{thm:Applications}.
Most of the results in this section are immediate consequences
of Theorem~\ref{thm:TreeGraded} in conjunction with results
from the literature.  Consequently, must of this section is expository in
nature.
The applications are divided into two subsections,
the first of which consists of
geometric results relating to quasi-isometries
and the second of which concerns consequences of relative hyperbolicity.

In Section~\ref{sec:Equivalent} we study the various equivalent
geometric formulations of isolated flats.
In particular we give the proof of Theorem~\ref{thm:GeometricEquivalence}.
We make use of these formulations
in the proof of Theorem~\ref{thm:TitsStructure},
which establishes the last remaining equivalence
(\ref{item:MainIsolated})~$\Longleftrightarrow$~(\ref{item:MainTits})
of Theorem~\ref{thm:MainEquivalence}.

\subsection{Acknowledgements}
\label{subsec:Thanks}

The first author would like to thank Benson Farb, Karen Vogtmann,
and Daniel Wise for many helpful conversations.
The first author is supported in part by an NSF Postdoctoral
Research Fellowship.
The second author is supported by NSF grant DMS-0204506.

\section{Background and preliminaries}
\label{sec:Background}

\subsection{Notation and terminology}
\label{subsec:Notation}

In this section, we establish notation and terminology that will be used
throughout this article.
We refer the reader to \cite{BH99} and \cite{Ballmann95}
for detailed introductions to $\CAT(0)$ spaces.
For an introduction to asymptotic cones, we refer the reader to
\cite{Gromov93} and \cite{KleinerLeeb97}.
The notion of a tree-graded space is due to \Drutu--Sapir
\cite{DrutuSapir05}.

Let $X$ be a $\CAT(0)$ space.  For points $p,q \in X$, we let
$[p,q]$ denote the unique geodesic connecting $p$ and~$q$.
If $c$ and $c'$ are geodesics emanating from a common point~$x$,
then $\angle_x(c,c')$ denotes the angle at~$x$ between $c$ and~$c'$.
If $x\ne y$ and $x\ne z$, then $\angle_x(y,z)$ denotes
the angle at~$x$ between $[x,y]$ and $[x,z]$.
The comparison angle at~$x$ corresponding to the angle $\angle_x(y,z)$
is denoted $\cangle_x(y,z)$.

For each $p \in X$, let $\Sigma_p X$ denote the space of directions
at~$p$ in~$X$ and let $\oa{px}$ denote the direction at~$p$
corresponding to the geodesic $[p,x]$.  The \emph{logarithm map}
\[
   \log_p \co X\setminus\{p\} \to \Sigma_p X
\]
sends each point~$x$ to the corresponding direction $\oa{px}$.

The ideal boundary of~$X$ equipped with the cone topology is denoted
$\boundary X$ and called the visual or geometric boundary.
The Tits angle metric on $\boundary X$ is denoted by $\tangle$
and the coresponding length metric by $\Td$.
The Tits boundary $\tits X$ is the ideal boundary of~$X$ equipped with
the metric $\Td$.

\begin{defn}\label{def:PointedHausdorff}
Let $X$ be a proper metric space.  The set of all closed subspaces
of~$X$ has a natural topology of \emph{Hausdorff convergence on
bounded sets} defined as follows.
For each closed set $C_0 \subseteq X$, each $x_0 \in X$, and
each positive $r$ and~$\epsilon$, define $U(C_0,x_0,r,\epsilon)$
to be the set of all closed subspaces $C \subseteq X$ such that
the Hausdorff distance between $C \cap B(x_0,r)$
and $C_0 \cap B(x_0,r)$ is less than~$\epsilon$.
The topology of Hausdorff convergence on bounded sets is the topology 
generated by the sets $U(\cdot,\cdot,\cdot,\cdot)$.
\end{defn}

\begin{defn}[Asymptotic cone]
Let $(X,d)$ be a metric space with a sequence of basepoints 
$(\star_n)$,
and let $(\lambda_n)$ be a sequence of numbers
called \emph{scaling constants} with $\ulim \lambda_n = \infty$.
The \emph{asymptotic cone} $\Cone_\omega(X,\star_n,\lambda_n)$
of $X$ with respect to $(\star_n)$ and $(\lambda_n)$
is the ultralimit $\ulim(X,\lambda_n^{-1} d,\star_n)$.
\end{defn}

Any ultralimit of a sequence of $\CAT(0)$
spaces is a complete $\CAT(0)$ space.
Any asymptotic cone of Euclidean space~$\E^k$ is isometric to~$\E^k$.
Let $X$ be any $\CAT(0)$ space with basepoints $(\star_n)$ and
scaling constants $(\lambda_n)$, and let $(F_n)$ be a sequence
of $k$-flats in~$X$ such that $\ulim \lambda_n^{-1} d(F_n,\star_n)$ is
finite.  Then the ultralimit of the sequence of embeddings
\[
   \bigl( F_n, \lambda_n^{-1} d_{F_n}, \pi_{F_n}(\star_n) \bigr)
   \inclusion  ( X, \lambda_n^{-1} d, \star_n )
\]
is a $k$-flat in
$\Cone_\omega(X,\star_n,\lambda_n)$.

\begin{defn}[Tree-graded]
Let $X$ be a complete geodesic metric space and let $\Set{P}$ be a collection
of closed geodesic subspaces of~$X$ called \emph{pieces}.
We say that $X$ is \emph{tree-graded with respect to~$\Set{P}$}
if the following two properties hold.
   \begin{enumerate}
   \item Every two distinct pieces have at most one common point.
   \item Every simple geodesic triangle in~$X$
         (a simple loop composed of three geodesics) lies inside one piece.
   \end{enumerate}
\end{defn}

\begin{defn}[Relatively hyperbolic spaces]
Let $X$ be a space and $\A$ a collection of subspaces of~$X$.
For each asymptotic cone $X_\om=\Cone_\om (X, \star_n, \lambda_n)$,
let $\A_\om$ be the collection of all subsets $A \subset X_\om$
of the form $A = \ulim A_n$ where $A_n \in \A$ and
$\ulim \lambda_n^{-1}\, d(A_n,\star_n) < \infty$.
Then $X$ is \emph{relatively hyperbolic with respect to~$\A$}
if, for every nonprincipal ultrafilter~$\omega$, each
asymptotic cone $X_\om$ is tree-graded with respect to~$A_\om$.
\end{defn}

\subsection{Triangles with small angular deficit}
\label{subsec:SmallDeficit}

In this subsection, we consider triangles whose angles are nearly
the same as the corresponding comparison angles.
We prove Proposition~\ref{prop:Subdivision},
which states that the family of such triangles is closed under 
certain elementary subdivisions, and Proposition~\ref{prop:SmallTriangle},
which provides a means of constructing numerous such triangles.

\begin{prop}\label{prop:Subdivision}
Let $\Delta(x,y,z)$ be a geodesic triangle in a $\CAT(0)$ space,
and choose $p \in [y,z] \setminus \{x,y,z\}$.
Suppose each angle of $\Delta(x,y,z)$
is within $\delta$ of the corresponding comparison angle.
Then each angle of $\Delta(x,y,p)$ and each angle of
$\Delta(x,p,z)$ is within $3\delta$ of the corresponding
comparison angle.

Furthermore, if each comparison angle of
$\Delta(x,y,z)$ is greater than $\theta$ for some $\theta>0$,
then $\angle_p(x,y)$, $\angle_p(x,z)$, and their comparison angles
lie in $(\theta,\pi-\theta)$.
\end{prop}

In order to prove the proposition, we first establish the following
lemma concerning subdivisions of a triangle with one angle close to its comparison angle.

\begin{lem}\label{lem:Subdivision}
Choose points $x,y,z \in X$ and $p\in [y,z] \setminus \{x,y,z\}$.
\begin{enumerate}
  \item\label{item:OuterCorner}
  Suppose $\angle_y(x,z)$ is within $\delta$ of the corresponding comparison
  angle. Then the quantities
  $\angle_y(x,p)$, \ $\cangle_y(x,p)$, and $\cangle_y(x,z)$
  are within~$\delta$.
  \item\label{item:SplitCorner}
  Suppose $\angle_x(y,z)$ is within $\delta$ of the corresponding comparison
  angle.  Then the same is true of both $\angle_x(y,p)$ and $\angle_x(p,z)$.
  Furthermore, the quantities
  $\cangle_x(y,z)$ and $\cangle_x(y,p) + \cangle_x(p,z)$
  are within~$\delta$.
\end{enumerate}
\end{lem}

\begin{proof}
The first assertion is immediate, since
\[
   \angle_y(x,z) = \angle_y(x,p)
     \le \cangle_y(x,p)
     \le \cangle_y(x,z)
     \le \angle_y(x,z) + \delta.
\]
For the second assertion, observe that
\begin{align*}
   \angle_x(y,z) &\le \angle_x(y,p) + \angle_x(p,z)   \\
                 &\le \cangle_x(y,p) + \cangle_x(p,z) \\
		 &\le \cangle_x(y,z)                  \\
		 &\le \angle_x(y,z) + \delta,
\end{align*}
where the inequality on the third line is due to Alexandrov's Lemma.
The assertion now follows from
the fact that each angle is no greater than its comparison angle.
\end{proof}

\begin{proof}[Proof of Proposition~\ref{prop:Subdivision}]
Lemma~\ref{lem:Subdivision} shows
that each of the angles
\[
   \angle_y(x,p), \quad \angle_z(x,p), \quad \angle_x(y,p),
   \quad \text{and} \quad \angle_x(p,z)
\]
is within~$\delta$ of its corresponding comparison angle.
But since $p \in [y,z]$, we also have
\begin{align*}
   \pi &\le \angle_p(x,y) + \angle_p(x,z) \\
       &\le \cangle_p(x,y) + \cangle_p(x,z) \\
       &= \bigl(\pi - \cangle_y(x,p) - \cangle_x(y,p) \bigr)
          + \bigl(\pi - \cangle_z(x,p) - \cangle_x(z,p) \bigr) \\
       &\le 2\pi + 3\delta - \cangle_y(x,z) - \cangle_x(y,z)
          - \cangle_z(x,y)  \\
       &= \pi + 3\delta,
\end{align*}
where the inequality on the fourth line is a consequence of
Lemma~\ref{lem:Subdivision}.
It now follows easily
that each of $\angle_p(x,y)$ and $\angle_p(x,z)$ is within $3\delta$
of its comparison angle.

To establish the last assertion, suppose each comparison
angle of $\Delta(x,y,z)$ is greater than~$\theta$.
Then
\[
   \angle_p(x,y) \le \cangle_p(x,y) \le \pi - \cangle_y(p,x)
     < \pi - \theta.
\]
Similarly, we see that $\angle_p(x,z) \le \cangle_p(x,z) < \pi - \theta$.
But $\angle_p(x,y) + \angle_p(x,z)$ is at least~$\pi$
by the triangle inequality for $\angle_p$.
Thus $\angle_p(x,y)$ and $\angle_p(x,z)$ are each greater than~$\theta$ as
desired.
\end{proof}

\begin{prop}\label{prop:SmallTriangle}
Let $c$ and $c'$ be geodesics with $c(0) = c'(0) = x$.
Then for sufficiently small $t$ and~$t'$,
each angle of $\Delta \bigl( x, c(t),c'(t') \bigr)$
is within $\delta$ of the corresponding comparison angle.
\end{prop}

The proof of Proposition~\ref{prop:SmallTriangle} uses the following lemma.

\begin{lemma}\label{lem:OneSideApproach}
Choose $x,y,z \in X$ and
let $c(t)$ be the geodesic $[x,y]$
parametrized so that $c(0)=x$.  Then
\[
   \lim_{t\to 0} \angle_{c(t)} (x,z)
    = \lim_{t\to 0} \cangle_{c(t)} (x,z)
    = \pi - \angle_x(y,z).
\]
\end{lemma}

\begin{proof}
Consider the comparison triangle $\bar\Delta\bigl(x,z,c(t)\bigr)$.
The comparison angle $\cangle_x \bigl( c(t),z \bigr)$ tends to
$\angle_x(y,z)$ as $t\to 0$
by the first variation formula \cite[II.3.5]{BH99}.
Since $\cangle_z \bigl( c(t),x \bigr)$ clearly tends to zero
and the angles
of $\bar\Delta$ sum to~$\pi$, we see that
\[
   \limsup_{t\to 0} \angle_{c(t)} (x,z) \le
   \lim_{t\to 0} \cangle_{c(t)} (x,z) = \pi - \angle_x(y,z).
\]
But we also have
\begin{align*}
   \liminf_{t\to 0} \angle_{c(t)} (x,z)
     &\ge \liminf_{t\to 0} \bigl( \pi - \angle_{c(t)} (y,z) \bigr) \\
     &=   \pi - \limsup_{t\to 0} \angle_{c(t)} (y,z) \\
     &\ge \pi - \angle_x (y,z),
\end{align*}
where the first inequality follows from the triangle inequality for
$\angle_{c(t)}$ and the last uses the upper semicontinuity of~$\angle$.
\end{proof}

\begin{proof}[Proof of Proposition~\ref{prop:SmallTriangle}]
By the definition of $\angle_x(c,c')$,
all points $y$ and~$z$ on $c$ and~$c'$ sufficiently close to~$x$
have the property that 
$\cangle_x(y,z)$ is within $\delta$ of $\angle_x(y,z)$.
Furthermore, by Lemma~\ref{lem:OneSideApproach} for each fixed $z$,
sliding $y$ toward~$x$ along~$c$, we can guarantee that
$\cangle_y(x,z)$ is within $\delta$ of $\angle_y(x,z)$.

Thus we may assume that the angles of $\Delta(x,y,z)$ at $x$ and~$y$
are within $\delta$ of the corresponding comparison angles.
Applying Lemma~\ref{lem:OneSideApproach} again, gives us a point
$z' \in [x,z]$ so that the angles of $\Delta(x,y,z')$ at $x$ and~$z'$
are within $\delta$ of their comparison angles.
If we consider the subdivision of $\Delta(x,y,z)$
into subtriangles $\Delta(x,y,z')$ and $\Delta(z',y,z)$
and apply Lemma~\ref{lem:Subdivision}(\ref{item:SplitCorner}),
it is clear that $\angle_y(x,z')$ is within $\delta$ of $\cangle_y(x,z')$
as well.
\end{proof}

\subsection{Periodic Euclidean subspaces and the Tits boundary}
\label{subsec:PeriodicEuclidean}

Let $Y$ be a proper $\CAT(0)$ space, and $F \subseteq Y$ a closed,
convex subset.  Let $G \subset \Isom(Y)$ be
any group of isometries of~$Y$.  We say that $F$ is \emph{$G$--periodic}
if $\Stab_G(F) \defeq \set{g \in G}{g(F)=F}$ acts cocompactly on~$F$.
We say that $F$ is \emph{periodic} if it is $\Isom(Y)$--periodic.

\begin{lemma}
\label{lem:HalfPlaneStickingUp}
Let $F$ be a periodic subspace of~$Y$ isometric to Euclidean space
$E^k$ for some $k \ge 1$.
Then either $\tits F$ is isolated in $\tits Y$,
or there is a flat half-plane $H\subset Y$ meeting $F$ orthogonally.
\end{lemma}

\begin{proof}
Assume $\tits F$ is not isolated in $\tits Y$, ie,
assume that there are sequences $\xi_k\in\tits Y\setminus\tits F$,
$\eta_k\in \tits F$ such that
\begin{equation}
\label{tanglera0}
\tangle(\xi_k,\eta_k)=\tangle(\xi_k,\tits F)\ra 0
\end{equation}
as $k\ra \infty$.    Pick $p\in F$.   Let $d_F$ denote the distance
function from $F$.  Then by (\ref{tanglera0})
there exists $c_k\ra 0$ such that for all
$x\in [p,\xi_k]$ we have $d(x,F)\leq c_k\, d(x,p)$.   Since
$d_F$ is convex, this also means that
when we restrict $d_F$ to the ray $[p,\xi_k]$ we get a function
whose left and right derivatives are $\leq c_k$.   Pick $R<\infty$,
and let $x_k\in [p,\xi_k]$ be the point on $[p,\xi_k]$ where the
distance to $F$ is $R$.   Now let $g_k\in \isom(Y)$ be a sequence of
isometries preserving $F$ such that $g_k(x_k)$ remains in a bounded subset
of $Y$, and pass to a subsequence so that $g_k(x_k)$
converges to some point $x_\infty$ and the rays $g_k\bigl([p,\xi_k]\bigr)$
converge to a complete geodesic  $\ga$ passing through $x_\infty$.
The derivative estimate on $d_F\restr{[p,\xi_k]}$ implies
that $d_F$ restricts to the constant function $R$ on $\ga$.
By the Flat Strip Theorem, we obtain a flat strip of width $R$
meeting $F$ orthogonally.  Since this works for any
$R<\infty$, we can use  the stabilizer of $F$ again to obtain a flat
half-plane leaving $F$ orthogonally.
\end{proof}

\subsection{Ultralimits of triangles}
\label{subsec:TriangleUltralimit}

This subsection is devoted to the proof of the following proposition.

\begin{prop}\label{prop:UlimEQ}
Let $X = \ulim(X_k,d_k,\star_k)$
be the ultralimit of a sequence of based $\CAT(0)$ spaces,
and let $x$, $y$, and~$z$ be points of~$X$ with
$x\ne y$ and $x\ne z$.
Let $(x_k)$, $(y_k)$, and $(z_k)$ be sequences representing
$x$, $y$, and~$z$ respectively.
Then there is a sequence $(x'_k)$ representing~$x$ such that
$x'_k \in [x_k,y_k]$ and such that, for all sequences
$(y'_k)$ and $(z'_k)$ representing $y$ and~$z$, we have
\[
   \ulim \angle_{x'_k} (y'_k,z'_k) = \angle_x(y,z).
\]
\end{prop}

Proposition~\ref{prop:UlimEQ} has the following immediate corollary
regarding ultralimits of geodesic triangles.

\begin{cor}\label{cor:TriangleUltralimit}
Let $X=\ulim(X_k,d_k,\star_k)$ be the ultralimit of a sequence
of $\CAT(0)$ spaces.
Let $(x_k)$, $(y_k)$, and $(z_k)$ be sequences representing three
distinct points $x$, $y$, and~$z$ in~$X$.
Then there are sequences $(x'_k)$, $(y'_k)$, and $(z'_k)$
also representing $x$, $y$, and~$z$ with the property that
$x'_k$ and $y'_k$ lie on the segment $[x_k,y_k]$
and each angle of the triangle $\Delta(x,y,z)$
is equal to the ultralimit of the corresponding angles of
$\Delta(x_k,y_k,z_k)$. \qed
\end{cor}

The proof of Proposition~\ref{prop:UlimEQ} uses the following two lemmas.

\begin{lemma}\label{lem:IndependentCorners}
Let $X=\ulim(X_k,d_k,\star_k)$ as above,
and choose points $x$, $y$, and~$z$ in~$X$
with $x\ne y$ and $x\ne z$.
If $(x_k)$, $(y_k)$, and $(z_k)$ represent $x$, $y$,
and~$z$ respectively, then the quantity
\[
   \ulim \angle_{x_k} (y_k,z_k)
\]
depends only on the choice of $(x_k)$ and not
on the choice of $(y_k)$ and $(z_k)$.
\end{lemma}

\begin{proof}
Fix a sequence $(x_k)$ representing~$x$, and
choose sequences $(y_k)$ and $(y'_k)$ representing~$y$.
Observe that, since $x\ne y$ we have
\[
   \ulim \frac{d_k(y_k,y'_k)}{d_k(y_k,x_k)} = 0.
\]
So the comparison angle $\cangle_{x_k} (y_k,y'_k)$ has ultralimit
zero, and therefore so does $\angle_{x_k} (y_k,y'_k)$.
The lemma now follows easily from the triangle inequality
for~$\angle_{x_k}$.
\end{proof}

\begin{lemma}\label{lem:UlimLE}
Suppose $X=\ulim(X_k,d_k,\star_k)$ as above,
and we have $(x_k)$, $(y_k)$, and
$(z_k)$ representing points $x$, $y$, and~$z$ in~$X$
with $x\ne y$ and $x\ne z$.
Then
\[
   \ulim \angle_{x_k}(y_k,z_k) \le \angle_x (y,z).
\]
\end{lemma}

\proof
Let $c_k$ and $c'_k$ be geodesic parametrizations of the segments
$[x_k,y_k]$ and $[x_k,z_k]$ such that $c_k(0)=c'_k(0)=x_k$.
Taking ultralimits of $c_k$ and $c'_k$ produces geodesic parametrizations
$c$ and~$c'$ of the segments $[x,y]$ and $[x,z]$.

For every positive~$t$, we have
\[
   \angle_{x_k}(y_k,z_k) = \angle_{x_k}\bigl( c_k(t),c'_k(t) \bigr)
   \le \cangle_{x_k} \bigl( c_k(t),c'_k(t) \bigr).
\]
Taking ultralimits of both sides of this inequality, we see that
\[
   \ulim \angle_{x_k}(y_k,z_k)
    \le \ulim \cangle_{x_k} \bigl( c_k(t),c'_k(t) \bigr)
    = \cangle_x \bigl( c(t),c'(t) \bigr).
\]
Since the previous inequality holds for all $t > 0$, it remains
true in the limit as $t\to 0$; ie,
$$
   \ulim \angle_{x_k} (y_k,z_k)
    \le \angle_x (y,z)  . \eqno{\qed}
$$

\begin{proof}[Proof of Proposition~\ref{prop:UlimEQ}]
Choose sequences $(x_k)$, $(y_k)$, and $(z_k)$ representing $x$, $y$,
and~$z$ respectively.  The choice of $(y_k)$ and $(z_k)$ is irrelevant
by Lemma~\ref{lem:IndependentCorners}.
Thus by Lemma~\ref{lem:UlimLE}, it suffices to find another
sequence $(x'_k)$ representing~$x$ such that
\begin{equation}\label{eqn:UlimGE}
   \ulim \angle_{x'_k} (y_k,z_k) \ge \angle_x(y,z).
\end{equation}

In $X$, choose points $w^n \in [x,y]$ with $w^n\to x$.
Since $[x,y]$ is the ultralimit of the segments $[x_k,y_k]$
we can choose points $w_k^n \in [x_k, y_k]$
so that for each~$n$ the sequence $(w_k^n)$ represents $w^n$.

The triangle inequality for $\angle_{w_k^n}$, together with
Lemmas \ref{lem:UlimLE}~and~\ref{lem:OneSideApproach} give that
\begin{align*}
   \ulim_k \angle_{w_k^n} (y_k,z_k)
     &\ge \ulim_k \bigl( \pi - \angle_{w_k^n} (x_k,z_k) \bigr) \\
     &\ge \pi - \angle_{w^n} (x,z) \\
     &= \angle_x (y,z) - \epsilon^{(n)},
\end{align*}
where $\epsilon^{(n)}$ is a constant depending on~$n$
such that $\epsilon^{(n)} \to 0$ as $n\to\infty$.
In other words, for each $k$ and $n$ there is a constant
$\epsilon_k^{(n)}$ with $\ulim_k \epsilon_k^{(n)} = \epsilon^{(n)}$
such that
\[
   \angle_{w_k^n} (y_k,z_k) \ge \angle_x(y,z) - \epsilon_k^{(n)}.
\]

The desired sequence $(x'_k)$ is constructed using a diagonal argument
as follows.
Set $\Set{S}_0 = \N$ and recursively define sets
$\Set{S}_n$ with $\omega(\Set{S}_n) = 1$ so that
$\Set{S}_n$ is contained in $\Set{S}_{n-1}$ but does not contain
the smallest element of $\Set{S}_{n-1}$, and so that
for each $k \in \Set{S}_n$ we have $\epsilon_k^{(n)} < 2\epsilon^{(n)}$.
Then $(\Set{S}_n)$ is a strictly decreasing sequence of sets
with empty intersection.

Set $x'_k = w_k^{n(k)}$, where $n(k)$ is the unique natural number such that
$k \in \Set{S}_{n(k)} \setminus \Set{S}_{n(k)+1}$.
By construction, $\ulim_k n(k) = \infty$.
Thus $(x'_k)$ represents~$x$, and
\[
   \ulim_k \epsilon_k^{(n(k))} \le \ulim_k 2\epsilon^{(n(k))} = 0,
\]
so $(x'_k)$ satisfies~(\ref{eqn:UlimGE}) as desired.
\end{proof}

\section{Isolated flats and tree-graded spaces}
\label{sec:IsolatedTreeGraded}

The main goal of this section is Theorem~\ref{thm:TreeGraded},
which establishes the equivalence
(\ref{item:MainIsolated})~$\Longleftrightarrow$~(\ref{item:MainCone})
of Theorem~\ref{thm:MainEquivalence}.

\subsection{Isolated flats}
\label{subsec:Isolated}

Throughout this subsection, let $X$ be a $\CAT(0)$ space and $\Gamma$
be a group acting geometrically on~$X$.
For the purposes of this section,
the term \emph{isolated flats} refers to the formulation (IF2).

\begin{lem}\label{lem:LocallyFinite}
If $X$ has isolated flats with respect to~$\F$, then $\F$ is locally
finite\textup{;}
in other words, only finitely many elements of~$\F$ intersect any given
compact set.
\end{lem}

\begin{proof}
It suffices to show that only finitely many $F \in \F$ intersect each
closed metric ball $\ol{B}(x,r)$.
Let $\F_0$ be the collection of all flats $F \in \F$ intersecting
this ball.
Choose $\kappa = \kappa(1)$ so that any intersection of $1$--neighborhoods
of distinct elements of~$\F$ has diameter less than~$\kappa$.
Any sequence of distinct elements in~$\F_0$ contains a subsequence $(F_i)$
that Hausdorff converges on bounded sets.
In particular, whenever $i$ and~$j$ are sufficiently large,
there are closed discs
$D_i \subset F_i$ and $D_j \subset F_j$ of radius~$\kappa$
whose Hausdorff distance is
less than~$1$, contradicting our choice of~$\kappa$.
\end{proof}

\begin{lem}\label{lem:PeriodicFlats}
Any locally finite, $\Gamma$--invariant collection~$\F$ of flats
in~$X$ has the following properties.
   \begin{enumerate}
   \item\label{item:FinOrbits}
   The elements of~$\F$ lie in only finitely many $\Gamma$--orbits,
   and the stabilizers of the $F \in \F$ lie in only finitely many
   conjugacy classes.
   \item\label{item:Periodic} Each $F \in \F$ is $\Gamma$--periodic.
   \end{enumerate}
\end{lem}

In particular, the lemma holds if $X$ has isolated flats with respect
to~$\F$.

\begin{proof}
Let $K$ be a compact set whose $\Gamma$--translates cover~$X$.
(\ref{item:FinOrbits}) follows easily from the fact that
only finitely many flats in~$\F$ intersect $K$.
Thus we only need to establish~(\ref{item:Periodic}).
For each $F \in \F$,
let $\{g_i\}$ be a minimal set of group elements such that the
translates $g_i(K)$ cover~$F$.
If the flats $g_i^{-1}(F)$ and $g_j^{-1}(F)$ coincide, then
$g_jg_i^{-1}$ lies in $H \defeq \Stab_{\Gamma}(F)$. 
It follows that the $g_i$ lie in only finitely many right cosets $Hg_i$.
In other words, the sets $g_i(K)$ lie in only finitely many $H$--orbits.
But any two $H$--orbits $Hg_i(K)$ and $Hg_j(K)$ lie at a finite Hausdorff
distance from each other.
Thus any $g_i(K)$ can be increased to a larger compact set $K'$ so that the
translates of $K'$ under~$H$ cover~$F$.
So $F$ is $\Gamma$--periodic, as desired.
\end{proof}

A flat in~$X$ is \emph{maximal} if it is not contained in a finite
tubular neighborhood of a higher dimensional flat.
It follows easily from (IF2)
that every flat in~$\F$ is maximal.

\begin{corollary}\label{cor:FlatStabs}
The set of stabilizers of flats $F \in \F$ is precisely the set~$\A$
of maximal virtually abelian subgroups of~$\Gamma$ with rank at least two.
\end{corollary}

\begin{proof}
By the Flat Torus Theorem, each subgroup $A \in \A$
stabilizes a flat.  Lemma~\ref{lem:PeriodicFlats}
and the Bieberbach Theorem show that the stabilizer of
each flat $F \in \F$ is virtually abelian.
The correspondence now follows easily from the maximality of
elements of~$\F$ and~$\A$.
\end{proof}

\subsection{Approximately Euclidean triangles}
\label{subsec:NearAFlat}

In this subsection we prove that, in a space with isolated flats,
a sufficiently large triangle whose vertex angles
are approximately equal to the corresponding comparison angles
must lie close to a flat.

\begin{lemma}\label{lem:DiscInFlat}
Let $X$ be any proper, cocompact metric space.
For each positive $r$ and~$\epsilon$, there is a constant $a=a(r,\epsilon)$
so that for each isometrically embedded $2$--dimensional
Euclidean disc~$C$ in~$X$
of radius~$a$,
the central subdisc of~$C$ of radius~$r$ lies in an $\epsilon$--tubular
neighborhood of a flat.
\end{lemma}

\begin{proof}
If not, then there would be constants $r$ and~$\epsilon$ and a sequence of
closed Euclidean $2$--dimensional discs
$C_a$ for $a=r$, $r+1$, $r+2,\dots$, such that $C_a$ has radius~$a$
and the central subdisc of $C_a$ of radius~$r$ does not lie in the
$\epsilon$--tubular neighborhood of any flat.
Applying elements of the cocompact isometry group and passing to a
subsequence, we may arrange that the discs $C_a$ Hausdorff converge
on bounded sets to a $2$--flat, contradicting our choice of $r$ 
and~$\epsilon$.
\end{proof}

\begin{lemma}
\label{lem:FlatSector}
Suppose $X$ has isolated flats.
There is a decreasing function $D_1=D_1(\theta)<\infty$ such that if
$S\subset X$ is a flat sector of angle $\theta>0$,
then $S\subset \nbd{F}{D_1(\theta)}$ for some $F\in\F$.
\end{lemma}

\begin{proof}
Pick positive constants $\th$ and~$r$, and set $\epsilon \defeq 1$.
Let $a=a(r,\epsilon)$ be the constant given by Lemma~\ref{lem:DiscInFlat},
and choose $\rho = \rho(\th,a)$ so that,
for any sector~$S$ with angle at least~$\th$,
the entire sector~$S$ lies inside a $\rho$--tubular neighborhood of
the subsector
\[
   S' \defeq \bigset{s\in S}{d(s,\boundary S) \ge a}.
\]
Note that for fixed $a$, the quantity $\rho(\theta,a)$
is a decreasing function of~$\theta$.
By Lemma~\ref{lem:DiscInFlat}, for each $s \in S'$ the intersection
of~$S$ with $B(s,r)$ is a flat disc contained in
the $\epsilon$--neighborhood of a flat.
By (IF2--\ref{item:AllFlats}),
this disc also lies in
the $(D+\epsilon)$--neighborhood of a flat $F \in \F$. 
Then (IF2--\ref{item:ControlledInt})
shows that $F$ will be independent of the choice of
$s\in S'$ when $r$ is sufficiently large.
Setting $D_1(\theta) \defeq D + \epsilon + \rho$ completes the proof.
\end{proof}

\begin{lemma}
\label{lem:CAngleFlat}
Let $X$ have isolated flats.
For all $\th_0>0$, $R<\infty$, there exist $\de_1=\de_1(\th_0,R)$,
$\rho_1=\rho_1(\th_0,R)$ such that if $p,x,y\in X$
satisfy $d(p,x),d(p,y)>\rho_1$ and
\[
  \th_0 < \angle_p(x,y)\leq \cangle_p(x,y) < \angle_p(x,y)+\de_1 < 
  \pi-\th_0,
\]
then there is a flat $F\in\F$ such that
\[
   \bigl([p,x]\cup [p,y]\bigr)\cap B(p,R)\subset \nbd{F}{D_1(\theta_0)}.
\]
\end{lemma}

\begin{proof}
If the lemma were false, we would have $\th_0>0$, $R<\infty$
and sequences $p_k,x_k,y_k\in X$ such that $d(p_k,x_k),d(p_k,y_k)\ra \infty$ and
\[
   \lim_{k\ra\infty}\angle_{p_k}(x_k,y_k)
     =\lim_{k\ra\infty}\cangle_{p_k}(x_k,y_k)
     \rdefeq\th\in[\th_0,\pi-\th_0]
\]
and such that there is no $F \in F$ with
\[
   \bigl([p_k,x_k]\cup [p_k,y_k]\bigr) \cap B(p_k,R)
     \subset \nbd{F}{D_1(\th_0)}.
\]
Applying the group~$\Gamma$ and passing to subsequences,
we may assume that there exist $p\in X$, $\xi,\eta\in\tits X$
such that $[p_k,x_k]\ra [p,\xi]$ and $[p_k,y_k]\ra [p,\eta]$.
By triangle comparison it follows that $\angle_p(\xi,\eta)=\tangle(\xi,\eta)=\th$,
so $[p,\xi]\cup [p,\eta]$ bounds a flat sector; then
Lemma~\ref{lem:FlatSector} gives a contradiction.
\end{proof}

The following lemma is an easy consequence of the $\CAT(0)$ inequality
and the Law of Cosines.
The proof is left as an exercise for the reader.

\begin{lem}[Convexity]\label{lem:Convexity}
Suppose $d(x,y)$ and $d(x,z)$ are at least~$D$
and $\angle_x(y,z)$ is at least $\theta$.
Then $d(y,z) \ge 2D \sin (\theta/2)$. \qed
\end{lem}

\begin{prop}\label{prop:TriangleNearFlat}
For all $\theta_0>0$ there are $\delta_2=\delta_2(\theta_0)>0$ and
$\rho_2=\rho_2(\theta_0)$ such that
if $x,y,z\in X$, all vertex angles and comparison angles
of $\Delta(x,y,z)$
lie in $(\theta_0,\pi-\theta_0)$, each vertex angle is within
$\delta_2$ of the corresponding comparison angle, and all
distances are greater than~$\rho_2$,
then
\[
   [x,y] \cup [x,z] \cup [y,z] \subset \nbd{F}{D_1(\theta_0)}
\]
for some flat $F \in \F$.
\end{prop}

\begin{proof}
Fix $\theta_0 > 0$ and let $R$ be sufficiently large that
each set of diameter at least $R/2$ lies in the
$D_1(\theta)$--tubular neighborhood of at most one flat $F \in \F$.
Let $\delta_1(\theta_0,R)$ and $\rho_1(\theta_0,R)$ be the constants
guaranteed by Lemma~\ref{lem:CAngleFlat},
and define
\[
   \rho_2 \defeq \max \left\{ 2\rho_1,
                      \frac{\rho_1}{\sin(\theta_0/2)} \right\}
     \quad \text{and} \quad
   \delta_2 \defeq \delta_1/3.
\]

Pick any triangle $\Delta(x,y,z)$
such that each vertex angle is within $\delta_2$ of the corresponding
comparison angle, all vertex angles and comparison angles
lie in $(\theta_0, \pi - \theta_0)$, and all side lengths are
greater than~$\rho_2$.
Choose an arbitrary point on one of the sides of~$\Delta$,
say $p \in [y,z]$.
Then $p$ divides $[y,z]$ into two segments, one of which, say $[p,y]$,
has length greater than $\rho_2/2$.

We now verify that the points $p$, $x$, and~$y$ satisfy the hypothesis
of Lemma~\ref{lem:CAngleFlat}.
First note that $\angle_p(x,y)$ and $\cangle_p(x,y)$ are within
$3\delta_2 = \delta_1$ of each other and both lie in
$(\theta_0,\pi-\theta_0)$ by Proposition~\ref{prop:Subdivision}.
By hypothesis
$d(p,y) > \rho_2/2 \ge \rho_1$.
Furthermore, observe that $d(p,y)$ and
$d(x,y)$ are both greater than $\rho_1/2\sin(\theta_0/2)$,
and that $\angle_p(x,y) \ge \theta_0$.
Thus by Lemma~\ref{lem:Convexity}, we have
$d(p,y) > \rho_1$.

Therefore by Lemma~\ref{lem:CAngleFlat} there is a flat $F_p \in \F$
such that
\[
   \bigl([p,x]\cup [p,y]\bigr)\cap B(p,R)\subset\bignbd{F_p}{D_1(\theta_0)}.
\]
But our choice of~$R$ guarantees that $F_p$
is independent of the choice of $p \in [x,y]\cup[x,z]\cup[y,z]$.
In other words, there is a single flat $F \in \F$ such that
\[
   [x,y]\cup[x,z]\cup[y,z] \subset \nbd{F}{D_1(\theta_0)},
\]
completing the proof.
\end{proof}

\subsection{Asymptotic cones are tree-graded}
\label{subsec:IsolatedIsTreeGraded}

Let $X_\om$ be an asymptotic cone $\Cone_\omega(X,\star_n,\lambda_n)$
where $X$ is $\CAT(0)$ with isolated flats. 
We let $\F_\omega$ denote the collection of flats $F$ in $X_\om$
of the form $F = \ulim F_n$ where $F_n \in \F$ and
$\ulim\lambda_n^{-1} \,d(F_n,\star) <\infty$.

\begin{lemma}
\label{lem:TriangleInFlat}
For all $\th_0>0$ there is a $\delta_3=\delta_3(\th_0)>0$ such that
if  $x,y,z\in X_\om$ are distinct, all vertex angles and comparison angles
of $\Delta(x,y,z)$ lie
in $(\th_0,\pi-\th_0)$, and each vertex angle is within
$\de_3$ of the corresponding comparison angle, then there
is a flat $F\in \F_\om$ containing $\{x,y,z\}$.  Furthermore,
if $F,F'\in \F_\om$ and $F\cap F'$ contains more than one
point, then $F=F'$.
\end{lemma}

\begin{proof}
Choose $\theta_0$, and let $\delta_2$ and $\rho_2$ 
be the constants provided by Proposition~\ref{prop:TriangleNearFlat}.
Set $\delta_3 := \delta_2$.
To prove the first assertion, choose $x,y,z \in X_\om$ as above
and apply Corollary~\ref{cor:TriangleUltralimit}
to get sequences $(x_k)$, $(y_k)$, and $(z_k)$ representing
$x$,~$y$, and~$z$ such that each vertex angle of $\Delta(x,y,z)$
is the ultralimit of the corresponding angle of $\Delta(x_k,y_k,z_k)$.
Then $\Delta(x_k,y_k,z_k)$ satisfies the hypothesis of
Proposition~\ref{prop:TriangleNearFlat} for $\omega$--almost all~$k$.
Consequently, $x$,~$y$, and~$z$ lie in a flat $F \in \F_\om$.

To prove the second assertion,
suppose flats $F,F' \in\F_\om$ have distinct points $x,y$ in common.
Let $(F_k)$ be a sequence in $\F$ with ultralimit~$F$,
and choose $z \in F'$ such that the triangle $\Delta(x,y,z)$ in~$F'$
is nondegenerate.
Let $(x_k)$, $(y_k)$, and $(z_k)$ be sequences in~$X$
representing $x$,~$y$, and~$z$, such that 
$x_k,y_k \in F_k$.
Since $[x_k,y_k]$ is contained in $F_k$ as well,
we may assume
by Corollary~\ref{cor:TriangleUltralimit} that
the angles of $\Delta(x,y,z)$ are ultralimits of the corresponding
angles of $\Delta(x_k,y_k,z_k)$.
Since $\Delta(x,y,z)$ is a flat triangle,
its vertex angles are equal to their comparison angles.
Therefore there exists $\theta_0 > 0$ such that
for $\om$--almost all~$k$ the triangle
$\Delta(x_k,y_k,z_k)$ satisfies
the hypothesis of Proposition~\ref{prop:TriangleNearFlat}.
Hence $\Delta(x_k,y_k,z_k)$
lies $D_1(\theta_0)$--close to a single flat $F_k''$.
As $[x_k,y_k] \subset F_k$,
we must have $F_k''=F_k$ for $\om$--almost all~$k$
by (IF2).
Hence $z\in F$.  Since a dense subset of $F'$ lies in $F$, it follows
that $F'\subseteq F$.
\end{proof}

\begin{prop}
\label{prop:DirectionComponents}
If $x\in X_\om$, then each connected component of\/
$\Si_{x} X_\om$ is either an isolated point or a sphere
of the form $\Si_{x} F$ for some flat $F \in \F_\om$ passing
through $x$.
The map $F \mapsto \Si_{x} F$ gives a one-to-one
correspondence between flats of $\F_\om$
passing through~$x$ and spherical components of $\Si_x F$.
Furthermore, if a direction $\oa{x y}$ lies in a spherical
component $\Si_{x} F$ then an initial segment of
$[x,y]$ lies in~$F$.
\end{prop}

\begin{proof}
Suppose $y,z\in X_\om$ and $0<\angle_{x}(y,z)<\pi$.
Then $\angle_x(y,z) \in (\theta, \pi - \theta)$ for some
positive~$\theta$.
Let $\delta_3 = \delta_3(\theta/8)$ be the constant given by
Lemma~\ref{lem:TriangleInFlat},
and let $\delta := \min \{\delta_3, \theta/4\}$.
Sliding $y$ and~$z$ toward~$x$ along the segments $[x,y]$ and $[x,z]$,
we may assume by Proposition~\ref{prop:SmallTriangle}
that the angles of $\Delta(x,y,z)$
are within~$\delta$ of their respective
comparison angles, and also that $d(x,y) = d(x,z)$.
Since $\delta \le \theta/4$, the angles of $\Delta(x,y,z)$
at $y$ and~$z$ lie in the interval $(\theta/8,\pi/2)$.
Lemma~\ref{lem:TriangleInFlat} now
implies that $x$,~$y$, and~$z$ lie in a common flat $F \in\F_\om$.
Thus the directions $\oa{x y}$ and $\oa{x z}$
lie in a sphere of the form $\Si_{x}F$ for some $F \in \F_\om$,
and each has an initial segment that lies in~$F$.

If $z'$ is any other point of~$X_\om$ with $0 < \angle_{x}(y,z') < \pi$,
then $[x,y]$ and $[x,z']$ have initial segments in a flat $F' \in \F_\om$.
It follows from Lemma~\ref{lem:TriangleInFlat} that $F=F'$.
Therefore $\Si_{x} F$ is a component of $\Si_{x}X_\om$.
Furthermore, every nontrivial component of $\Si_x X_\om$ arises in this
manner, and $\Si_x F = \Si_x F'$ implies $F=F'$.
\end{proof}

\begin{corollary}\label{cor:ConstantProjection}
Let $\pi_{F} \co X_\om \to F$
be the nearest point projection of $X_\om$ onto a flat $F \in \F_\om$.
Then $\pi_F$ is locally constant on $X_\om \setminus F$.
In other words, $\pi_F$ is constant on each component of
$X_\om \setminus F$.
\end{corollary}

\begin{proof}
Choose $s \in X_\om \setminus F$ and let $x = \pi_F(s)$.
Then $\angle_x (s,F)$ is at least $\pi/2$.
So $\log_x(s)$ lies in a different component of $\Sigma_x X_\om$
from the tangent sphere $\Sigma_x F$
by Proposition~\ref{prop:DirectionComponents}.
By continuity of $\log_x$, if $U$ is any connected neighborhood of~$s$
not containing~$x$, the image $\log_x(U)$ is disjoint
from the component $\Sigma_x F$.
Thus each point of $\log_x(U)$ is at an angular distance~$\pi$
from $\Sigma_x F$.
Hence for each $s'\in U$, we have $\pi_{F}(s') = x$.
\end{proof}

\begin{lemma}\label{lem:CutPoint}
If $p$ lies on the interior of the segment
$[x,y]\subset X_\om$, then $x$ and $y$ lie
in distinct components of $X_\om \setminus \{p\}$
unless $p$ is contained in an open subarc of $[x,y]$ that lies in
a flat $F\in\F_\om$.
\end{lemma}

\begin{proof}
If the segments $[p,x]$ and $[p,y]$ do not have initial segments
in a common flat,
then the directions $\oa{px}$ and $\oa{py}$
lie in distinct components of $\Sigma_p X_\om$.
The result follows immediately from the continuity of $\log_p$.
\end{proof}

\begin{lemma}
\label{lem:LoopInFlat}
Every embedded loop in $X_\om$ lies in some flat $F\in \F_\om$.
\end{lemma}

\begin{proof}
Let $\gamma$ be an embedded loop containing points $x \ne y$.
By Lemma~\ref{lem:CutPoint} and the fact that $\gamma$ has
no cut points, it follows that the geodesic $[x,y]$
lies in some flat $F \in \F_\om$.
Let $\beta$ be a maximal open subpath of~$\gamma$ in the complement
of~$F$.  It follows from Corollary~\ref{cor:ConstantProjection} that
$\beta$ projects to a constant under $\pi_{F}$.
Hence its endpoints coincide, which is absurd.
\end{proof}

\begin{thm}\label{thm:TreeGraded}
Let $X$ be a $\CAT(0)$ space admitting a proper, cocompact,
isometric action of a group~$\Gamma$.
If $X$ has isolated flats with respect to a family of flats~$\F$,
then $X$ is relatively hyperbolic with respect to~$\F$.
\end{thm}

\begin{proof}
Each flat $F \in \F_\om$ is a closed convex subspace of $X_\om$.
By Lemma~\ref{lem:TriangleInFlat}, any two flats
$F \ne F' \in \F_\om$ intersect in at most one point.
Furthermore, Lemma~\ref{lem:LoopInFlat} shows that every embedded
geodesic triangle in $X_\om$ lies in a flat $F \in \F_\om$.
\end{proof}

\section{Applications}
\label{sec:Applications}

Theorem~\ref{thm:TreeGraded} has many immediate consequences,
which we examine in the present section.
In particular, we show the equivalence
(\ref{item:MainCone})~$\Longleftrightarrow$~(\ref{item:MainRelHyp})
of Theorem~\ref{thm:MainEquivalence} and prove the various
parts of Theorem~\ref{thm:Applications}.
These applications are divided into two subsections, the first of which
deals with geometric properties of spaces with isolated flats.
The second subsection is concerned with relative hyperbolicity
and its ramifications.

\subsection{Geometric invariants of spaces with isolated flats}
\label{subsec:GeometricInvariants}

The following theorem is an immediate consequence of
Theorem~\ref{thm:TreeGraded} together with
\cite[Proposition~5.4]{DrutuSapir05}.

\begin{thm}[Quasiflat]
Let $X$ be a $\CAT(0)$ space with isolated flats.
For each constants $L$ and~$C$, there is a constant $M = M(L,C)$
so that every quasi-isometrically embedded flat in~$X$ lies
in an $M$--tubular neighborhood of some flat $F \in \F$. \qed
\end{thm}

We will now begin to explore the relationship between
the geometry of $X$ and the algebra of a group $\Gamma$
acting geometrically on~$X$.

\begin{prop}\label{prop:TGFlatAbelian}
Suppose $\Gamma$ acts geometrically on a $\CAT(0)$ space $X$.
The following are equivalent.
   \begin{enumerate}
   \item \label{item:TGFlats}
   $X$ is relatively hyperbolic with respect to a $\Gamma$--invariant
   collection~$\F$ of flats.
   \item \label{item:TGAbelian}
   $\Gamma$ is a relatively hyperbolic metric space with respect to the
   left cosets of a collection~$\A$ of virtually abelian subgroups
   of rank at least two.
   \end{enumerate}
\end{prop}

\begin{proof}
The action of $\Gamma$ on~$X$ induces a quasi-isometry $\Gamma\to X$.
Since being relatively hyperbolic with respect to quasiflats
is a geometric property
\cite[Theorem~5.1]{DrutuSapir05},
it suffices to show that this quasi-isometry
takes $\F$ to~$\A$ if we assume either (\ref{item:TGFlats})
or (\ref{item:TGAbelian}).

By Theorem~\ref{thm:TreeGraded}, if $X$ satisfies (\ref{item:TGFlats})
then $X$ has isolated flats with respect to~$\F$.
So by Lemma~\ref{lem:PeriodicFlats} and the Bieberbach Theorem,
each $F\in \F$ is $\Gamma$--periodic with virtually abelian stabilizer.
Let $\A$ be a collection of representatives of the 
conjugacy classes of stabilizers of flats $F\in\F$.
Then the quasi-isometry $\Gamma\to X$
sends the elements of~$\F$ to the left cosets of elements of~$\A$. 

Conversely, suppose $X$ satisfies (\ref{item:TGAbelian}).
Let $\F_0$ be a collection of flats stabilized by the elements of~$\A$,
as guaranteed by the Flat Torus Theorem.
Then $\Gamma\to X$ again sends the $\Gamma$--translates of elements
of~$\F_0$ to the left cosets of elements of~$\A$.
\end{proof}

Recall that being virtually abelian is a quasi-isometry invariant
of finitely generated groups.
Therefore, by \cite[Corollary~5.19]{DrutuSapir05},
satisfying Property~(\ref{item:TGAbelian}) of
Proposition~\ref{prop:TGFlatAbelian} is also a quasi-isometry invariant
of finitely generated groups.
We therefore have the following corollary.

\begin{cor}
Let $X_1$ and $X_2$ be quasi-isometric $\CAT(0)$ spaces
admitting geometric actions by groups $\Gamma_1$ and $\Gamma_2$
respectively.
If $X_1$ has isolated flats, then $X_2$ also has isolated flats. \qed
\end{cor}

\begin{defn}[Fellow travelling relative to flats.]
\label{def:RelativeFellowTravelling}
A pair of paths
\[
   \alpha \co [0,a]\to X \quad \text{and}
   \quad \alpha' \co [0,a']\to X
\]
in
a $\CAT(0)$ space \emph{$L$--fellow travel relative to a sequence of flats}
$(F_1,\dots,F_n)$
if there are partitions
\[
  0 = t_0 \le s_0 \le t_1 \le s_1 \le \dots \le t_n \le s_n = a
\]
and
\[
  0 = t'_0\le s'_0\le t'_1\le s'_1\le \dots \le t'_n\le s'_n= a'
\]
so that for $0\le i \le n$ the Hausdorff distance between the
sets $\alpha \bigl( [t_i,s_i] \bigr)$ and
$\alpha' \bigl( [t'_i,s'_i] \bigr)$ is at most~$L$,
while for $1\le i \le n$ the sets
$\alpha \bigl( [s_{i-1},t_i] \bigr)$
and $\alpha' \bigl( [s'_{i-1},t'_i] \bigr)$ lie in an $L$--neighborhood of
the flat $F_i$.

We will frequently say that paths \emph{$L$--fellow travel relative
to flats} if they $L$--fellow travel relative to some sequence of flats.
\end{defn}

\begin{defn}\label{def:RelFTP}
A $\CAT(0)$ space $X$ satisfies the
\emph{Relative Fellow Traveller Property}
if for each choice of constants $\lambda$ and~$\epsilon$
there is a constant $L = L(\lambda,\epsilon,X)$ such that
$(\lambda,\epsilon)$--quasigeodesics in~$X$
with common endpoints $L$--fellow travel relative to flats.
\end{defn}

The following is an immediate consequence of 
Theorem~\ref{thm:TreeGraded} together with
the Morse Property for relatively hyperbolic spaces
proved by \Drutu--Sapir in
\cite[Theorem~4.25]{DrutuSapir05}.

\begin{prop}\label{prop:RFTP}
If $X$ has isolated flats, then it also has the Relative
Fellow Traveller Property. \qed
\end{prop}

The Relative Fellow Traveller Property was previously established by
Epstein for the truncated hyperbolic space associated to
a finite volume cusped hyperbolic manifold
\cite[Theorem~11.3.1]{ECHLPT92}.
In the context of $2$--dimensional $\CAT(0)$ spaces,
Hruska showed that isolated flats is equivalent to the
Relative Fellow Traveller Property and also to the Relatively
Thin Triangle Property, which states that each geodesic triangle
is thin relative to a flat in a suitable sense \cite{Hruska2ComplexIFP}.

In \cite{HruskaGeometric}, Hruska proved several results
about $\CAT(0)$ spaces with isolated flats under the additional assumption
that the Relative Fellow Traveller Property holds.
With Proposition~\ref{prop:RFTP}, we can now drop the Relative Fellow Traveller Property as a hypothesis in each of those theorems.
In the remainder of this subsection, we list some immediate consequences
of Proposition~\ref{prop:RFTP} together with \cite{HruskaGeometric}.

A subspace $Y$ of a geodesic space~$X$ is \emph{quasiconvex}
if there is a constant~$\kappa$ so that every geodesic in~$X$
connecting two points of~$Y$ lies in the $\kappa$--tubular
neighborhood of~$Y$.
Let $\rho \co G \to \Isom(X)$ be a geometric action of a group on a
$\CAT(0)$ space.
A subgroup $H \le \Gamma$ is \emph{quasiconvex with respect to~$\rho$}
if any (equivalently ``every'') $H$--orbit is a quasiconvex subspace
of~$X$.

In the word hyperbolic setting, Short has shown that
any finitely generated $H\le \Gamma$ is quasiconvex
if and only if it is \emph{undistorted} in~$\Gamma$;
in other words, the inclusion $H \inclusion \Gamma$ is a quasi-isometric
embedding \cite{Short91}.
It is easy to see that this equivalence does not extend
to the general $\CAT(0)$ setting (see for instance \cite{HruskaGeometric}).
However, in the presence of isolated flats, we have the following theorem.

\begin{theorem}\label{thm:QC-CAT0}
If $\rho$ is a geometric action
of $\Gamma$ on a $\CAT(0)$ space~$X$ with isolated flats,
then a finitely generated subgroup
$H \le \Gamma$ is quasiconvex with respect to~$\rho$ if and only
if $H$ is undistorted in~$\Gamma$. \qed
\end{theorem}

We remark that Osin and \Drutu--Sapir have examined related phenomena
involving undistorted subgroups and relatively quasiconvex subgroups
of relatively hyperbolic groups in
\cite[\S 4.2]{OsinRelHyp} and \cite[\S 8.3]{DrutuSapir05}.
In general, it is not true that all undistorted subgroups of
a relatively hyperbolic group are quasiconvex.  One can only conclude that
they are ``relatively quasiconvex'' with respect to the parabolic subgroups.
In two slightly different contexts,
Osin and \Drutu--Sapir study undistorted subgroups that have finite
intersections with all parabolic subgroups, essentially
proving that such subgroups are quasiconvex in the standard sense
and are word hyperbolic.

The main difference between Theorem~\ref{thm:QC-CAT0} and the work
of Osin and \Drutu--Sapir is that Theorem~\ref{thm:QC-CAT0}
applies to undistorted subgroups with arbitrary intersections with the
parabolic subgroups.  The $\CAT(0)$ geometry inside the flats plays a
crucial role in proving quasiconvexity in the isolated flats setting
in \cite{HruskaGeometric}.

Now suppose $\Gamma$ acts geometrically on two $\CAT(0)$ spaces
$X_1$ and~$X_2$.
In the word hyperbolic setting, Gromov observed that the geometric
boundary is a group invariant, in the sense that $\boundary X_1$ and 
$\boundary X_2$ are $\Gamma$--equivariantly homeomorphic.
More generally, a quasiconvex subgroup of~$\Gamma$ is again
word hyperbolic, and its boundary embeds equivariantly into the boundary
of~$\Gamma$ \cite{Gromov87}.

Croke--Kleiner have shown that the homeomorphism type
of the geometric boundary is not a group invariant in the general
$\CAT(0)$ setting \cite{CrokeKleiner00}.
In fact Wilson has shown that the Croke--Kleiner construction
produces a continuous family of homeomorphic $2$--complexes whose universal 
covers have nonhomeomorphic geometric boundaries \cite{Wilson05}.

In the presence of isolated flats, the situation is quite
similar to the hyperbolic setting.
However, a subtle complication arises from the fact that
it is currently unknown whether a quasiconvex subgroup of a
$\CAT(0)$ group is itself $\CAT(0)$.  (A group is $\CAT(0)$ if it admits
a geometric action on a $\CAT(0)$ space.)

\begin{theorem}[Boundary of a quasiconvex subgroup is well-defined]
Let $\rho_1$ and $\rho_2$ be geometric actions of groups
$\Gamma_1$ and $\Gamma_2$
on $\CAT(0)$ spaces $X_1$ and $X_2$ both with isolated flats.
For each $i$, let $H_i \le \Gamma_i$ be a subgroup quasiconvex with respect
to~$\rho_i$.
Then any isomorphism  $H_1 \to H_2$ induces an equivariant
homeomorphism $\Lambda H_1 \to \Lambda H_2$
between the corresponding limit sets. \qed
\end{theorem}

Setting $H_i \defeq \Gamma_i$ gives the following corollary.

\begin{cor}[Boundary is well-defined]
If $\Gamma$ acts geometrically
on a $\CAT(0)$ space~$X$ with isolated flats, then
the geometric boundary of~$X$ is a group invariant of~$\Gamma$. \qed
\end{cor}

\subsection{Relative hyperbolicity and its consequences}
\label{subsec:RelativeHyperbolicity}

The notion of a relatively hyperbolic group was first proposed by
Gromov in \cite{Gromov87}.  A substantially different definition was
proposed and studied by Farb in \cite{Farb98}.  Gromov's approach was
later formalized by Bowditch and shown to be equivalent to Farb's
definition in \cite{BowditchRelHyp}.  For our purposes,
the most useful characterization of relatively hyperbolic groups is given
by the following theorem.  (The converse implication is proved by
\Drutu--Sapir in \cite{DrutuSapir05}.  The direct implication
is proved by Osin--Sapir in an appendix to \cite{DrutuSapir05}.)

\begin{theorem}[\Drutu--Osin--Sapir]
\label{thm:TG=RelHyp}
A finitely generated group~$\Gamma$ is relatively hyperbolic
with respect to a collection~$\A$ of finitely generated subgroups
if and only if
$\Gamma$ is relatively hyperbolic as a metric space with respect to the set
of left cosets of the subgroups $A \in \A$.
\end{theorem}

The subgroups $A \in \A$ are called \emph{peripheral subgroups}
of the relatively hyperbolic structure.

The equivalence
(\ref{item:MainCone})~$\Longleftrightarrow$~(\ref{item:MainRelHyp})
of Theorem~\ref{thm:MainEquivalence}
follows immediately from Theorem~\ref{thm:TG=RelHyp} and 
Proposition~\ref{prop:TGFlatAbelian}.

A group $\Gamma$ satisfies the Tits Alternative if every subgroup
of~$\Gamma$
either is virtually solvable or contains a nonabelian free subgroup.
For $\CAT(0)$ groups this property is equivalent to the Strong Tits Alternative,
which states that every subgroup either is virtually abelian or contains
a nonabelian free subgroup.  The Strong Tits Alternative was proved
for word hyperbolic groups by Gromov \cite[3.1.A]{Gromov87},
and has also been established for $\CAT(0)$ groups acting geometrically
on certain real analytic $4$--manifolds \cite{Xie4Manifold04} or on
cubical complexes
(proved for special cases in \cite{BS99} and \cite{XieTitsAltSquared04},
and in full generality in \cite{SageevWiseTitsAlt}).
However, it is still unknown whether
the Strong Tits Alternative holds for arbitrary $\CAT(0)$ groups.

Theorem~\ref{thm:MainEquivalence}, together with work of Gromov
and Bowditch, allows us to extend the Strong Tits Alternative to the
class of $\CAT(0)$ groups acting on spaces with isolated flats.
More precisely,
Gromov shows in \cite[8.2.F]{Gromov87} that any properly discontinuous
group of isometries of a proper $\delta$--hyperbolic space
either is virtually cyclic, is parabolic (ie,
fixes a unique point at infinity),
or contains a nonabelian free subgroup.
According to Bowditch (elaborating on an idea due to Gromov),
each relatively hyperbolic group acts
properly discontinuously on a proper $\delta$--hyperbolic
space such that the maximal parabolic subgroups
are exactly the conjugates of the peripheral subgroups \cite{BowditchRelHyp}.
The following result follows immediately.

\begin{thm}[Gromov--Bowditch]
Let $\Gamma$ be relatively hyperbolic with respect to a collection of
subgroups that each satisfy the \textup{[}Strong\textup{]} Tits Alternative.
Then $\Gamma$ also satisfies the \textup{[}Strong\textup{]}
Tits Alternative.
\end{thm}

In the isolated flats
setting, the peripheral subgroups are virtually abelian.  We therefore
conclude the following.

\begin{thm}[Isolated flats $\implies$ Tits Alternative]
A group~$\Gamma$ that acts geometrically
on a $\CAT(0)$ space with isolated flats
satisfies the Strong Tits Alternative.  In other words,
every subgroup $H \le \Gamma$ either is virtually abelian or contains a
nonabelian free subgroup. \qed
\end{thm}

All word hyperbolic groups are biautomatic \cite{ECHLPT92}.
As with the Tits Alternative, it is currently unknown whether this result
also holds for arbitrary $\CAT(0)$ groups.
In fact, the only nonhyperbolic spaces where biautomaticity 
was previously known
are complexes built out of very restricted shapes of cells.
For instance,
Gersten--Short established biautomaticity for $\CAT(0)$ groups
acting geometrically on $\CAT(0)$ $2$--complexes of type
$\tilde{A}_1 \times \tilde{A}_1$, $\tilde{A}_2$, $\tilde{B}_2$, and 
$\tilde{G}_2$
in \cite{GerstenShort90Automatic,GerstenShort91Automatic}.
In particular, Gersten--Short's work includes
$\CAT(0)$ square complexes and
$2$--dimensional Euclidean buildings.
Niblo--Reeves proved biautomaticity for groups
acting geometrically on arbitrary
$\CAT(0)$ cube complexes in \cite{NibloReeves98}.

Rebbechi shows in \cite{Rebbechi01} that a group which is hyperbolic
relative to a collection of biautomatic subgroups is itself biautomatic.
Since finitely generated, virtually abelian groups are biautomatic, we
obtain the following immediate corollary to 
Theorem~\ref{thm:MainEquivalence}.

\begin{thm}[Isolated flats $\implies$ biautomatic]
If $\Gamma$ acts geometrically
on a $\CAT(0)$ space with isolated flats, then $\Gamma$ is biautomatic. \qed
\end{thm}

Many $\CAT(0)$ $2$--complexes with isolated flats have
irregularly shaped cells.  Such $2$--complexes provide new examples
of biautomatic groups.  The following is an elementary
construction of such a $2$--complex.

\begin{example}[Irregularly shaped cells]
Let $X$ be a compact hyperbolic surface and $\gamma$ a geodesic loop in~$X$
representing a primitive conjugacy class.  Let $T$
be any flat $2$--torus containing a simple geodesic loop~$\gamma'$.
Form a $2$--complex~$Y$ by gluing $X$ and~$T$ along the curve
$\gamma = \gamma'$.  Then the universal cover of~$Y$ has isolated flats.
Typically $\gamma \subset X$ will intersect itself many times, subdividing
$X$ into a large number of irregularly shaped pieces.  
\end{example}

\section{Equivalent formulations of isolated flats}
\label{sec:Equivalent}

In this section we prove Theorem~\ref{thm:GeometricEquivalence},
which establishes the equivalence of the various geometric formulations
of isolated flats discussed in the introduction.
We then use Theorem~\ref{thm:GeometricEquivalence}
to prove the remaining equivalence
(\ref{item:MainIsolated})~$\Longleftrightarrow$~(\ref{item:MainTits})
of Theorem~\ref{thm:MainEquivalence},
which characterizes spaces with isolated flats in terms of
their Tits boundaries.

\subsection{Proof of Theorem~\ref{thm:GeometricEquivalence}}
\label{subsec:FirstProperties}

The implications (\ref{item:UniformlyThin})~$\Longrightarrow$
(\ref{item:Thin}) $\Longrightarrow$ (\ref{item:Slim}) are immediate
from the definitions.  We prove (\ref{item:IF1})
$\Longrightarrow$ (\ref{item:Thin}) in Proposition~\ref{prop:isolated=>thin}.
We then show in Proposition~\ref{prop:slim=>isolated}
that (\ref{item:Slim}) implies both (\ref{item:IF1})
and~(\ref{item:UniformlyThin}).
Finally we prove (\ref{item:Slim}) $\Longrightarrow$ (\ref{item:IF2})
and (\ref{item:IF2}) $\Longrightarrow$ (\ref{item:IF1})
in Propositions \ref{prop:slim=>controlled} 
and~\ref{prop:controlled=>isolated} respectively.

\begin{lemma}\label{lem:LocallyFiniteII}
Any closed, isolated $\Gamma$--invariant subset $\F \subseteq \flat(X)$
is locally finite.
\end{lemma}

\begin{proof}
If $K\subset X$ is compact, then
\[
   \set{F\in\F}{F\cap \bar K\neq\emptyset}
\]
is a closed subset of the compact set
\[
  \set{F\in\flat(X)}{F\cap \bar K\neq\emptyset},
\]
and is therefore compact; but all its elements are isolated,
so it is also finite.   Thus $\F$ is a locally finite collection.
\end{proof}

The following is an immediate corollary of Lemmas
\ref{lem:LocallyFiniteII} and~\ref{lem:PeriodicFlats}.

\begin{cor}\label{cor:PeriodicII}
Let $\F$ be a closed, isolated $\Gamma$--invariant subset of $\flat(X)$.
Then the elements of~$\F$ lie in only finitely many $\Gamma$--orbits
and each $F \in \F$ is $\Gamma$--periodic.\qed
\end{cor}

\begin{lemma}\label{lem:Maximal}
Suppose $X$ has \textup{(IF1)} with respect to the family~$\F$.
We may assume without loss of generality
that each element of~$\F$ is maximal.
\end{lemma}

\begin{proof}
By Corollary~\ref{cor:PeriodicII}, the elements of~$\F$ lie in only
finitely many $\Gamma$--orbits.
Thus we can delete $\Gamma$--orbits of flats from~$\F$
until it becomes a minimal $\Gamma$--invariant subset
satisfying (IF1).
But each element of a minimal~$\F$ is clearly maximal.
\end{proof}

\begin{prop}\label{prop:isolated=>thin}
If $X$ has \textup{(IF1)} then $X$ has thin parallel sets.
\end{prop}

\begin{proof}
It is sufficient to
consider parallel sets of geodesics lying in flats $F\in \F$.
Let $\Xi$ be the collection of pairs $(F,S)$ where $F\in\F$,
$S\subset F$ is isometric to $\E^k$ for some $k\geq 1$,
and the parallel set of $S$ is not contained in a finite neighborhood
of~$F$.   Suppose $\Xi$ is nonempty, and choose an element $(F,S)\in\Xi$
which is ``maximal'' in the sense that if $(F',S')\in\Xi$
and $\dim F'\geq \dim F$ and $\dim S'\geq \dim S$, then
$\dim F'=\dim F$ and $\dim S'=\dim S$.

\emph{Step~1.}  We first show that $S=F$.
Assume by way of contradiction
that $\dim S<\dim F$, so that $\dim\bar F\geq 1$.

We have a Euclidean product decomposition $\P(S)=S\times Y$ where
$Y\subset X$ is convex.
Let $\pi \co \P(S)\to Y$ be the projection onto the second factor,
and define
$\bar F\subset Y = \pi(F)$.  Then we have $F=S\times \bar F$.
Note that $Y$ is not contained in a finite neighborhood of $\bar F$.
So by applying the cocompact stabilizer of~$F$ and a convergence
argument, we may assume without loss of generality that
$\tits Y\setminus\tits \bar F\neq\emptyset$.

Suppose $\tits \bar F$ is not a component of $\tits Y$.  Then
we may apply Lemma~\ref{lem:HalfPlaneStickingUp}  to see that there
is a flat half-plane $H\subset Y$ which meets $\bar F$ orthogonally.
But then $S'\defeq S\times \D H\subset F$ is a subflat of $F$ with
$\dim S'>\dim S$ and $\P(S')\supset S\times H$ is not contained in a finite
neighborhood of $F$.  This contradicts the choice of the pair
$(F,S)$, and so $\tits \bar F$ must be a component of $\tits Y$.

Pick $\xi\in \tits Y\setminus\tits \bar F$.  Then there is a geodesic
$\gamma\subset Y$ with $\tits\gamma=\{\xi,\eta\}$
where $\eta\in\tits \bar F$.
By (IF1) there is a flat $F'\in\F$
such that $S\times \ga\subset \nbd{F'}{D}$;
let $\bar F'$ be the projection of $F'$ to $Y$.  Let $\al\subset
\bar F$ be a geodesic with $\eta\in\tits\al$, and pick $x\in F$
such that $\bar x\defeq \pi(x)\in \al$.  By the Bieberbach Theorem,
the stabilizer of~$F$ contains a subgroup~$A$ which acts
cocompactly by translations on $F$;
therefore we can find an
infinite sequence $g_k\in A$ such that $d\bigl(g_k(\al),\al\bigr)$ is
uniformly bounded, and $g_k\bar x$ converges
to $\eta_-\in \geo X$, where $\tits\al=\{\eta,\eta_-\}$.   All the flats
$g_k(F')$ pass through some ball around $x$, so the collection
$\bigl\{g_k(F')\bigr\}$ is finite.  It follows that for each $\eps>0$,
there is a translation
$g\in A$ which preserves $F'$,  and which translates $\bar F$
in a direction $\eps$--close to $\eta_-$; in particular $g$ does not
preserve $S$.  If $\tau\subset \bar F$ is a geodesic translated
by $g$, then $S'\defeq S\times\tau$ is a Euclidean subset with
$\dim S'>\dim S$, and $\P(S')\supset F'$ is not contained in a finite neighborhood
of $F$.  This contradicts the choice of $(F,S)$.  Therefore we
have $S=F$.

\emph{Step~2.}  We now know that $S=F$.  So $\P(F)=F\times Y$
where $Y$ is unbounded.
Pick $p\in F$, and let $[p,x_k] \subset \P(F)$ be a sequence of geodesics
in the $Y$--factor through~$p$ such that
$d(x_k,F)\to\infty$.
Let $y_k$ be the midpoint of the
segment $[p,x_k]$, and let $F_k\subset \P(F)$ be the flat
parallel to $F$ passing through $y_k$.
By Lemma~\ref{lem:Maximal}, $F$ is maximal.
Thus the parallel flat $F_k$ is maximal as well.
Applying isometries $g_k \in \Gamma$ and passing to a subsequence,
we can arrange that the maximal flats $g_k(F_k)$ converge to a flat $\bar{F}$
that is not maximal.

Each maximal flat $g_k(F_k)$ is at a Hausdorff distance at most~$D$
from some $\hat{F}_k \in \F$ with $\dim\hat{F}_k = \dim{F}$.
Since $g_k(F_k)$ converges
and $\F$ is locally finite, we can pass to a further subsequence so that
the flats $\hat{F}_k$ are equal to a single flat $\hat{F} \in \F$.
Thus the Hausdorff distance between $\hat{F}$ and~$\bar{F}$
is at most~$D$.
But $\hat{F}$ is maximal and $\bar{F}$ is not, which is a contradiction.
\end{proof}

\begin{prop}\label{prop:slim=>isolated}
If $X$ has slim parallel sets,
then $X$ has uniformly thin parallel sets
and satisfies \textup{(IF1)}.
\end{prop}

The proposition is an immediate consequence of the following lemma
in the case $n=2$.

\begin{lemma}\label{lem:slim=>isolated}
Suppose $X$ has slim parallel sets.
Then for each $n\ge 2$ there is a closed, isolated $\Gamma$--invariant set
$\F_n \subset \flat(X)$ and a constant~$D_n$ with the following properties.
\begin{enumerate}
   \item Each element of~$\F_n$ is a maximal flat of dimension at least~$n$.
   \item Each $k$--flat in~$X$ with $k\ge n$
   is contained in the $D_n$--tubular neighborhood of some $F \in \F_n$.
   \item \label{item:UniformDn}
   If $\gamma\subseteq F$ is a geodesic contained in a flat
   $F \in \F_n$
   then the parallel set of~$\gamma$ lies in the $D_n$--neighborhood
   of~$F$.
\end{enumerate}
\end{lemma}

\begin{proof}
Since there is an upper bound on the dimensions of all flats in~$X$,
the lemma is satisfied for sufficiently large~$n$
by the set $\F_n \defeq \emptyset$.
We will prove the lemma by decreasing induction on~$n$.
Suppose the lemma is true for $n+1$.

If $F$ is a maximal $n$--flat, then by hypothesis its parallel set has the form
$F \times K_F$ for some compact set $K_F$.
Define $\F_n$ to be the union of
$\F_{n+1}$ and the set of all $n$--flats of the form $F \times \{c_F\}$,
where $F$ is maximal and $c_F$ is the circumcenter of~$K_F$.
Note that $\F_n$ is $\Gamma$--invariant.

\emph{Step 1.}  We will first show that the set of all maximal flats
of dimension $n$ is closed in $\flat(X)$.
If not, then there is a sequence
$(F_i)$ of maximal $n$--flats converging to a nonmaximal $n$--flat~$F$,
which must lie in a finite tubular neighborhood of a higher dimensional
flat $\hat{F} \in \F_{n+1}$.
It follows from (\ref{item:UniformDn}) that
$F \subset \nbd{\hat{F}}{D_{n+1}}$.
Since $F$ is not maximal, it is not parallel to any~$F_i$.
So $F_i$ is not contained in any finite tubular neighborhood of~$\hat{F}$.
Consequently, for any $r>D_{n+1}$,
if we set
\[
   C_i \defeq F_i \cap \nbd{\hat{F}}{r},
\]
then $C_i$ is a convex, open, proper subset of~$F_i$ whose inscribed
radius~$\rho_i$ tends to infinity as $i \to \infty$.
Let $p_i \in F_i$ be a point on the boundary of $C_i$
that is also on the boundary of an inscribed $n$--dimensional Euclidean disc
in~$C_i$ of radius~$\rho_i$.
By Corollary~\ref{cor:PeriodicII}, we know that $\hat{F}$ has a cocompact
stabilizer.
So we may apply elements $g_i \in \Stab(\hat{F})$ and pass to
a subsequence so that the $g_i(p_i)$ converge to a point $p$,
the $g_i(F_i)$ converge to a flat~$\bar{F}$ containing~$p$,
the $g_i(C_i)$ converge to an open halfspace $H \subset \bar{F}$
that lies in the closed $r$--tubular neighborhood of~$\Hat{F}$,
and each point of the bounding hyperplane $\boundary H$ lies
at a distance exactly~$r$ from $\hat{F}$.
Therefore $\boundary H$ lies in the parallel set of a geodesic in~$\hat{F}$
and $r>D_{n+1}$,
which contradicts the inductive hypothesis applied to~$\hat{F}$.

\emph{Step~2.} Our next goal is to show that a sequence of
pairwise nonparallel
maximal flats cannot converge in $\flat(X)$.
Suppose by way of contradiction that $(F_i)$ is such a sequence
that converges in $\flat(X)$ to~$F$.
By Step~1, we know that $F$ is maximal.
Fix a positive~$r$, and
define $C_i$ and $p_i$ in $F_i$ as in Step~1.
As $F$ is $n$--dimensional, it is not yet known to be periodic.
Nevertheless, we may apply isometries $g_i \in \Gamma$
and pass to a subsequence to arrange that $g_i(p_i)$
converges to a point~$p$,
that $g_i(F_i)$ converges to a flat $\bar{F}$,
that $g_i(F)$ converges to a flat $\hat{F}$,
and that $g_i(\boundary C_i)$ converges to an $(n-1)$--dimensional
subflat~$S$ of~$\bar{F}$ each point of which is at a distance exactly~$r$
from~$\hat{F}$.
By Step~1, the limiting flats $\bar{F}$ and~$\hat{F}$ are again maximal.
But $\bar{F} \cup \hat{F}$ is contained in the parallel set of
any geodesic line $\gamma \subset S$.
It follows from slim parallel sets that
$\tits \bar{F} = \tits \hat{F}$, and hence that $\bar{F}$ and $\hat{F}$
are parallel and separated by a distance exactly~$r$.

Since $r$ is arbitrary, the preceding paragraph produces for each
positive~$r$
a pair of parallel maximal $n$--flats $\bar{F}_r$ and~$\hat{F}_r$
at a Hausdorff distance exactly~$r$.
The convex hull of $\bar{F}_r \cup \hat{F}_r$ is isometric to
a product $\bar{F}_r \times [-r/2,r/2]$.
Applying isometries in~$\Gamma$ again, we can arrange that
as $r \to \infty$
the maximal flats $\bar{F}_r \times \{0\}$ converge to an $n$--flat
which is not maximal, contradicting Step~1.

Thus $\F_n$ is a closed, isolated, $\Gamma$--invariant set in $\flat(X)$
containing exactly one flat from each parallel class
of maximal flats of dimension $\ge n$.
By Corollary~\ref{cor:PeriodicII},
the elements of $\F_n$ lie in only finitely
many $\Gamma$--orbits and are each $\Gamma$--periodic.
Another convergence argument using the cocompact stabilizer
of each $F \in \F_n$ can now be used to uniformly bound the thickness of
the parallel sets $\P(\gamma)$ for each geodesic $\gamma \subset F$.
The lemma follows immediately, since
each $k$--flat with $k\ge n$ lies inside one of these parallel sets.
\end{proof}

\begin{prop}\label{prop:slim=>controlled}
Let $X$ have slim parallel sets.  Then $X$ has \textup{(IF2)}.
\end{prop}

\begin{proof}
Consider the set $\F \defeq \F_2$ given by Lemma~\ref{lem:slim=>isolated}.
Then $\F$ is locally finite, consists of only maximal flats, 
and also has the property that no two flats in~$\F$ are parallel.
Pick $\rho<\infty$ and a sequence
of pairs $F_k\ne F_k'$ in~$\F$ such that
$\diam \bigl( \nbd{F_k}{\rho} \cap \nbd{F_k'}{\rho} \bigr)\ra\infty$.
After passing to subsequences and
applying a sequence of isometries, we may assume that
$F_k\to F_\infty\in\F$, $F_k'\to F'_\infty\in\F$, and
$\nbd{F_\infty}{\rho}\cap \nbd{F'_\infty}{\rho}$
contains a complete geodesic~$\ga$.
Since $\F$ is locally finite, we in fact
have $F_k=F_\infty$, $F_k'=F'_\infty$ for all sufficiently large~$k$.

But $F_\infty \cup F'_\infty$ lies in the parallel set of~$\gamma$.
As in the proof of Lemma~\ref{lem:slim=>isolated}, it follows that
$F_\infty$ and $F'_\infty$ are parallel, contradicting the fact that
$F_k \ne F'_k$.
\end{proof}

\begin{prop}\label{prop:controlled=>isolated}
If $X$ has \textup{(IF2)} with respect to~$\F$,
then $X$ also has \textup{(IF1)} with respect to~$\F$.
\end{prop}

\begin{proof}
The only nontrivial fact to check is that
a sequence of distinct flats in $\F$ cannot converge in $\flat(X)$.
But this follows immediately from Lemma~\ref{lem:LocallyFinite}.
\end{proof}

\subsection{The structure of the Tits boundary}
\label{subsec:TitsBoundaryStructure}

In this section we show that $X$ has isolated flats if and only if
the components of its Tits boundary are isolated points and standard
spheres.

\begin{prop}
\label{prop:RaysFlat}
If $X$ has isolated flats, then
for each $\th_0>0$ there is a positive constant $\de_4=\de_4(\th_0)$
such that whenever
$p\in X$ and $\xi,\eta\in\tits X$ satisfy
\begin{equation}
\label{eqn:RayAnglePinch}
   \theta_0 \leq \angle_p(\xi,\eta)
      \leq \tangle(\xi,\eta)
      \leq \angle_p(\xi,\eta) + \delta_4
      \leq \pi - \theta_0,
\end{equation}
then there is a flat $F \in \F$ so that
\[
  [p,\xi]\cup [p,\eta]\subset \bignbd{F}{D_1(\th_0)}.
\]
\end{prop}

\begin{proof}
Choose $\theta_0$ positive, and let $R$ be sufficiently large that
any set of diameter at least $R/4$ is contained in the
$D_1(\theta_0)$--tubular neighborhood of at most one flat in~$\F$.
Let $\delta_1 = \delta_1(\theta_0,R)$ and
$\rho_1 = \rho_1(\theta_0,R)$ be the constants given by
Lemma~\ref{lem:CAngleFlat}, and set $\delta_4:=\delta_1$.
Now suppose $p\in X$ and $\xi,\eta\in\tits X$
satisfy
\begin{equation}
\label{eqn:RayAnglePinch1}
   \theta_0 \leq \angle_p(\xi,\eta)
      \leq \tangle(\xi,\eta)
      \leq \angle_p(\xi,\eta) + \delta_1
      \leq \pi - \theta_0.
\end{equation}
Then we have
\[
   \angle_p(\xi,\eta)
     \le \cangle_p \bigl( c(\rho_1),c'(\rho_1) \bigr)
     \le \tangle (\xi,\eta),
\]
where $c$ and $c'$ are geodesic parametrizations of $[p,\xi]$
and $[p,\eta]$ respectively.
So by Lemma~\ref{lem:CAngleFlat} there is a flat $F_p \in \F$ so that
\[
   B(p,R)\cap\bigl([p,\xi]\cup [p,\eta]\bigr)\subset \bignbd{F}{D_1(\th_0)}.
\]

Notice that
(\ref{eqn:RayAnglePinch1}) remains valid if we replace $p$ with
any point $x \in [p,\xi]\cup [p,\eta]$.
Thus we can apply the preceding argument to the rays $[x,\xi]$ and $[x,\eta]$
to produce for each $x$ a flat $F_x \in \F$ such that
\[
   B(x,R) \cap \bigl( [x,\xi] \cup [x,\eta] \bigr)
     \subset \bignbd{F_x}{D_1(\theta_0)}.
\]
But our choice of~$R$ guarantees that all the flats $F_x$
are the same.  So $[p,\xi] \cup [p,\eta]$ lies in the
$D_1(\theta)$--tubular neighborhood of a single flat $F \in \F$.
\end{proof}

The following proposition was
established by Schroeder \cite[Appendix~4.E]{BGS85}
in the Riemannian setting and extended to
arbitrary $\CAT(0)$ spaces by Leeb \cite{Leeb98}.
In the proposition, a \emph{standard sphere} is a sphere isometric to a
unit sphere in Euclidean space.  A \emph{standard hemisphere}
is isometric to a hemisphere of a standard sphere.

\begin{prop}[Schroeder, Leeb]
\label{prop:StandardSphere}
Let $X$ be a proper $\CAT(0)$ space, and $S=S^{k-1}$ an isometrically
embedded standard sphere in $\tits X$.
Then one of the following holds.
   \begin{enumerate}
   \item There is a $k$--flat $F$ in~$X$ such that $\tits F = S$.
   \item There is an isometric embedding of a standard
   hemisphere $H \subset S^k$
   in $\tits X$ such that the topological boundary of~$H$ maps
   to $S$.
   \end{enumerate}
\end{prop}

\begin{thm}\label{thm:TitsStructure}
Let $X$ be any $\CAT(0)$ space admitting a geometric group action.
Then $X$ has isolated flats if and only if
each connected component of $\tits X$ either is an isolated point
or is isometric to a standard sphere.
Furthermore, the spherical components of $\tits X$
are precisely the Tits boundaries of the flats $F \in \F$.
\end{thm}

\begin{proof}
Suppose $X$ has isolated flats.
By (IF2), it is clear that if $F,F' \in \F$
are distinct then $\tits F \cap \tits F' = \emptyset$.
If $\xi,\eta\in\tits X$ and $0<\tangle(\xi,\eta)<\pi$,
then we can find $\th_0>0$ and $p\in X$ such that
(\ref{eqn:RayAnglePinch}) holds.
Hence $[p,\xi]\cup [p,\eta]\subset \nbd{F}{D_1(\theta_0)}$
for some $F\in\F$, which means that $\{\xi,\eta\}\subset\tits F$.
Thus $\tits X$ has the desired structure.

Now suppose every component of $\tits X$ is either an isolated point
or a standard sphere.
We will show that $X$ has slim parallel sets.
By Proposition~\ref{prop:StandardSphere},
any spherical component of $\tits X$ is the boundary of a maximal flat.
It follows that the Tits boundary of any maximal flat~$F$
is a component of $\tits X$.  But for each geodesic $\gamma \subset F$,
the Tits boundary of $\P(\gamma)$ is a connected set
containing $\tits F$.  Hence $\tits \P(\gamma) = \tits F$.
\end{proof}

\appendix
\section{Appendix: Isolated subspaces}
\label{sec:Appendix}
\setcounter{equation}{0}

\vspace{-3mm}
\textsl{By Mohamad Hindawi, G~Christopher Hruska and Bruce Kleiner}
\vspace{4mm}

Hindawi proved in \cite{Hindawi05}
(independent of the main body of this article)
that the universal cover of
a closed, real analytic, nonpositively curved $4$--manifolds whose Tits
boundary has no nonstandard components exhibits similar behavior to
$\CAT(0)$ spaces with isolated flats.

In this appendix, we show that many of the results of this article
extend to a more general class of $\CAT(0)$ spaces, including
the manifolds studied by Hindawi.
The main idea is to weaken the hypothesis $\textup{(IF2)}$ so that
the elements of~$\F$, instead of being flats, are closed convex subspaces.
This general idea was introduced as a generalization
of isolated flats by Kapovich--Leeb in
\cite{KapovichLeeb95}.

\begin{theorem}\label{thm:AppMain}
Let $X$ be a $\CAT(0)$ space, and $\Ga\acts X$ be
a geometric action.
Suppose $\F$ is a
$\Ga$--invariant collection of closed convex subsets.
Assume that
\begin{enumerate}
\renewcommand{\theenumi}{\textup{\Alph{enumi}}}
\item \label{item:AppAllFlats}
There is a constant $D<\infty$ such that
each flat $F\subseteq X$ lies in a $D$--tubular neighborhood of
some $C\in \F$.
\item \label{item:AppControlledInt}
For each positive $r<\infty$ there is a constant $\rho=\rho(r)<\infty$
so that for any two distinct elements $C, C' \in \F$ we have
\[
   \diam \bigl( \nbd{C}{r} \cap \nbd{C'}{r} \bigr) < \rho.
\]
\end{enumerate}
Then\textup{:}
\begin{enumerate}
\item \label{item:AppPeriodic}
The collection $\F$ is locally finite,
there are only finitely many $\Ga$--orbits in~$\F$,
and each $C \in \F$ is $\Ga$--periodic.

\item \label{item:AppTitsBoundary}
Every connected component of $\tits X$ containing more
than one point is contained in $\tits C$, for some $C\in \F$.

\item \label{item:AppConeDirections}
For every asymptotic cone $X_\om$ of $X$, and every $p\in X_\om$,
each connected component of $\Si_p X_\om$ containing more
than one point is contained in $\Si_p C_\om$, for
some $C_\om\in F_\om$.

\item \label{item:AppTreeGraded}
Every asymptotic cone $X_\om$ is tree-graded with respect
to the collection $\F_\om$.

\item \label{item:AppRelHyp}
$\Ga$ is hyperbolic relative to the collection of
stabilizers of elements of $\F$.  

\item \label{item:AppBoundaryWellDefined}
Suppose the stabilizers of elements of $C\in \F$
are $\CAT(0)$ groups with very well-defined boundary
\textup{(}see below\textup{)}.
Then $\Ga$ has a very well-defined boundary.
\end{enumerate}
\end{theorem}

\begin{defn}
A group $\Gamma$ has a \emph{very well-defined boundary}
if, whenever $\Gamma$ acts geometrically on a pair of $\CAT(0)$ spaces
$X$ and $Y$, any $\Gamma$--equivariant quasi-isometry $X \to Y$
sends each geodesic to within a uniformly bounded Hausdorff distance
of a geodesic.
\end{defn}

Examples of groups with very well-defined boundary include
word hyperbolic groups, virtually abelian groups,
and any group acting geometrically on a product
$X_1 \times \dots \times X_k$ where each $X_i$ is either an irreducible
higher rank symmetric space or an irreducible higher rank Euclidean building.

\begin{proof}[Proof of Theorem~\ref{thm:AppMain}]
The point is that
many of the arguments in this article do not use the assumption
that the elements of~$\F$ are flats in an essential way.

The proof of~(\ref{item:AppPeriodic}) is identical to that of Lemmas
\ref{lem:LocallyFinite} and~\ref{lem:PeriodicFlats}.
The proof of~(\ref{item:AppTitsBoundary})
follows from the arguments in Sections \ref{subsec:NearAFlat}
and~\ref{subsec:TitsBoundaryStructure}.
Assertions (\ref{item:AppConeDirections}), (\ref{item:AppTreeGraded}),
and (\ref{item:AppRelHyp})
follow using arguments similar to those in Sections
\ref{sec:IsolatedTreeGraded} and~\ref{sec:Applications}.

Assertion (\ref{item:AppBoundaryWellDefined}) follows using arguments
similar to those in \cite{HruskaGeometric}.
Note that the arguments in \cite{HruskaGeometric} use the
Relative Fellow Traveller Property as a hypothesis.
However, this property is a special case
of \Drutu--Sapir's Morse Property for
relatively hyperbolic spaces
\cite{DrutuSapir05}.
The arguments of \cite{HruskaGeometric} go through
essentially unchanged if one uses the Morse Property
in place of the Relative Fellow Traveller Property
throughout.
\end{proof}

A partial converse to the preceding theorem is given by the following
result due to \Drutu--Sapir \cite{DrutuSapir05}.

\begin{thm}
Suppose $\Gamma$ is relatively hyperbolic with respect to
a finite collection of subgroups $\{P_i\}$.
Fix a finite generating set for~$\Gamma$, and consider the corresponding
word metric on $\Gamma$.
Then each $P_i$ is quasiconvex in~$\Gamma$ \textup{(}with
respect to the word metric\textup{)}.
Furthermore,
\begin{enumerate}
\item every quasiflat in~$\Gamma$ lies in a uniformly bounded
tubular neighborhood of some coset $gP_i$, and
\item for each positive $r < \infty$, there is a constant
$\rho=\rho(r) < \infty$ so that for any two distinct cosets $g_1P_1$
and $g_2 P_2$, we have
\[
   \diam \bigl( \nbd{g_1P_1}{r} \cap \nbd{g_2P_2}{r} \bigr) < \rho.
\]
\end{enumerate}
\end{thm}

\def\polhk#1{\setbox0=\hbox{#1}{\ooalign{\hidewidth
  \lower1ex\hbox{$\,\lhook$}\hidewidth\crcr\unhbox0}}}
  \def\RomanianComma#1{\setbox0=\hbox{#1}{\ooalign{\hidewidth
  \lower1.2ex\hbox{$\mspace{1mu}^{,}$}\hidewidth\crcr\unhbox0}}}


\begin{thebibliography}{}

\bibitem{AlexanderBergBishop93}
\textbf{S\,B Alexander}, \textbf{I\,D Berg}, \textbf{R\,L Bishop},
  \emph{Geometric curvature bounds in {R}iemannian manifolds with boundary},
  Trans. Amer. Math. Soc. 339 (1993) 703--716 \MR{MR1113693}

\bibitem{AlibegovicBestvinaLimit}
\textbf{E Alibegovi{\'{c}}}, \textbf{M Bestvina}, \emph{Limit groups are
  {$\CAT(0)$}}, e-print (2004) \arxiv{math.GR/0410198}

\bibitem{Ballmann95}
\textbf{W Ballmann}, \emph{Lectures on spaces of nonpositive curvature},
  volume~25 of \emph{DMV Seminar}, with an appendix by M~Brin,
  Birkh\"auser Verlag, Basel (1995)
  \MR{1377265}
  

\bibitem{BallmannBrin94}
\textbf{W Ballmann}, \textbf{M Brin}, \emph{Polygonal complexes and
  combinatorial group theory}, Geom. Dedicata 50 (1994) 165--191
  \MR{1279883}

\bibitem{BGS85}
\textbf{W Ballmann}, \textbf{M Gromov}, \textbf{V Schroeder}, \emph{Manifolds
  of nonpositive curvature}, Birkh\"auser, Boston, Mass. (1985)
  \MR{0823981}

\bibitem{BS99}
\textbf{W Ballmann}, \textbf{J {\'S}wi{\polhk{a}}tkowski}, \emph{On groups
  acting on nonpositively curved cubical complexes}, Enseign. Math. (2) 45
  (1999) 51--81
  \MR{1703363}

\bibitem{Benakli94}
\textbf{N Benakli}, \emph{Polygonal complexes. {I}. {C}ombinatorial and
  geometric properties}, J. Pure Appl. Algebra 97 (1994) 247--263\relax
  \MR{1314577}
  
\bibitem{Bowditch93}
\textbf{B\,H Bowditch}, \emph{Geometrical finiteness for hyperbolic groups}, J.
  Funct.\ Anal. 113 (1993) 245--317\relax \MR{94e:57016}

\bibitem{BowditchRelHyp}
\textbf{B\,H Bowditch}, \emph{Relatively hyperbolic groups}, preprint,
  University of Southampton (1999)

\bibitem{Bridson95}\textbf{M\,R Bridson}, \emph{On the existence of
flat planes in spaces of nonpositive curvature}, Proc.\ Amer.\ Math.\ Soc. 123
(1995) 223--235 \MR{1273477}

\bibitem{BH99}
\textbf{M\,R Bridson}, \textbf{A Haefliger}, \emph{Metric spaces of
  non-positive curvature}, Springer-Verlag, Berlin (1999) \MR{1744486}

\bibitem{Buyalo98}
\textbf{S\,V Buyalo}, \emph{Geodesics in {H}adamard spaces}, Algebra i Analiz
  10 (1998) 93--123; {E}nglish translation in
  {St. Petersburg Math. J.} 10 (1999) 293--313 \MR{MR1629391}

\bibitem{ChatterjiRuane05}
\textbf{I Chatterji}, \textbf{K Ruane}, \emph{Some geometric groups with
  {R}apid {D}ecay}, Geom. Funct. Anal. 15 (2005) 311--339

\bibitem{CrokeKleiner00}
\textbf{C\,B Croke}, \textbf{B Kleiner}, \emph{Spaces with nonpositive
  curvature and their ideal boundaries}, Topology 39 (2000) 549--556
  \MR{1746908}

\bibitem{CrokeKleiner02}
\textbf{C\,B Croke}, \textbf{B Kleiner}, \emph{The geodesic flow of a
  nonpositively curved graph manifold}, Geom. Funct. Anal. 12 (2002)
  479--545
  \MR{1924370}

\bibitem{DrutuSapirRD05}
\textbf{C Dru{\RomanianComma{t}}u}, \textbf{M Sapir}, \emph{Relatively
  hyperbolic groups with {R}apid {D}ecay property}, Int. Math. Res. Not.
  (2005) 1181--1194

\bibitem{DrutuSapir05}
\textbf{C Dru{\RomanianComma{t}}u}, \textbf{M Sapir}, \emph{Tree-graded spaces
  and asymptotic cones of groups}, with an appendix
  by D~Osin and M~Sapir, Topology 44 (2005) 959--1058

\bibitem{Eberlein73}
\textbf{P Eberlein}, \emph{Geodesic flows on negatively curved manifolds.
  {I}{I}}, Trans.\ Amer.\ Math.\ Soc. 178 (1973) 57--82
  \MR{0314084}

\bibitem{ECHLPT92}
\textbf{D\,B\,A Epstein}, \textbf{J\,W Cannon}, \textbf{D\,F Holt},
  \textbf{S\,V\,F Levy}, \textbf{M\,S Paterson}, \textbf{W\,P Thurston},
  \emph{Word processing in groups}, Jones and Bartlett Publishers, Boston, MA
  (1992)
  \MR{1161694}

\bibitem{Farb98}
\textbf{B Farb}, \emph{Relatively hyperbolic groups}, Geom.\ Funct.\ Anal. 8
  (1998) 810--840
  \MR{1650094}

\bibitem{GerstenShort90Automatic}
\textbf{S\,M Gersten}, \textbf{H\,B Short}, \emph{Small cancellation theory and
  automatic groups}, Invent. Math. 102 (1990) 305--334
  \MR{1074477}

\bibitem{GerstenShort91Automatic}
\textbf{S\,M Gersten}, \textbf{H\,B Short}, \emph{Small cancellation theory and
  automatic groups. {II}}, Invent. Math. 105 (1991) 641--662
  \MR{1117155}

\bibitem{Gromov81}
\textbf{M Gromov}, \emph{Hyperbolic manifolds, groups and actions}, from:
  ``Riemann surfaces and related topics\textup{:} Proceedings of the 1978 Stony
  Brook Conference \textup{(}State Univ. New York, Stony Brook, N.Y.,
  1978\textup{)}'', (I Kra, B Maskit, editors), Ann. of Math. Stud. 97,
  Princeton Univ. Press, Princeton, N.J. (1981)  183--213\relax \MR{MR624814}

\bibitem{Gromov87}
\textbf{M Gromov}, \emph{Hyperbolic groups}, from: ``Essays in group theory'',
  (S\,M Gersten, editor), Springer, New York (1987)  75--263
  \MR{0919829}

\bibitem{Gromov93}
\textbf{M Gromov}, \emph{Asymptotic invariants of infinite groups}, from:
  ``Geometric group theory, Vol.~2 \textup{(}Sussex, 1991\textup{)}'', (G\,A
  Niblo, M\,A Roller, editors), Cambridge Univ.\ Press, Cambridge (1993)
  1--295
  \MR{1253544}

\bibitem{GrovesHopfian}
\textbf{D Groves}, \emph{Limits of \textup{(}certain\textup{)} {$\CAT(0)$}
  groups, {II}\textup{:} {T}he {H}opf property and the shortening
  argument}, \arxiv{math.GR/0408080}

\bibitem{Haglund91}
\textbf{F Haglund}, \emph{Les poly\`edres de {G}romov}, C. R. Acad. Sci. Paris
  S\'er. I Math. 313 (1991) 603--606
  \MR{1133493}

\bibitem{Heber87} \textbf{J Heber}, \emph{Hyperbolische geodatische {R}aume},
Diplomarbeit, Univ.\ Bonn (1987)

\bibitem{Hindawi05}
\textbf{M\,A Hindawi}, \emph{Large scale geometry of $4$--dimensional closed
  nonpositively curved real analytic manifolds}, Int. Math. Res. Not.  (2005)
  1803--1815

\bibitem{Hruska2ComplexIFP}
\textbf{G\,C Hruska}, \emph{Nonpositively curved $2$--complexes with isolated
  flats}, Geom. Topol. 8 (2004) 205--275 \MR{2033482}

\bibitem{HruskaGeometric}
\textbf{G\,C Hruska}, \emph{Geometric invariants of spaces with isolated
  flats}, Topology 44 (2005) 441--458 \MR{MR2114956}

\bibitem{KapovichLeeb95}
\textbf{M Kapovich}, \textbf{B Leeb}, \emph{On asymptotic cones and
  quasi-isometry classes of fundamental groups of $3$--manifolds}, Geom.\
  Funct.\ Anal. 5 (1995) 582--603
  \MR{1339818}

\bibitem{KleinerLeeb97}
\textbf{B Kleiner}, \textbf{B Leeb}, \emph{Rigidity of quasi-isometries for
  symmetric spaces and {E}uclidean buildings}, Inst. Hautes \'Etudes Sci. Publ.
  Math. 86 (1997) 115--197
  \MR{1608566}

\bibitem{Leeb98}
\textbf{B Leeb}, \emph{A characterization of irreducible symmetric spaces and
  Euclidean buildings of higher rank by their asymptotic geometry},
  Habilitation, Univ. Bonn (1998)
  \MR{1934160}

\bibitem{Moussong88}
\textbf{G Moussong}, \emph{Hyperbolic {C}oxeter groups}, PhD thesis, Ohio
  State Univ. (1988)

\bibitem{NibloReeves98}
\textbf{G\,A Niblo}, \textbf{L\,D Reeves}, \emph{The geometry of cube complexes
  and the complexity of their fundamental groups}, Topology 37 (1998)
  621--633
  \MR{1604899}

\bibitem{OsinRelHyp}
\textbf{D\,V Osin}, \emph{Relatively hyperbolic groups\textup{:} {I}ntrinsic
  geometry, algebraic properties, and algorithmic
  problems}, to appear in Mem.\ Amer.\ Math.\ Soc. \arxiv{math.GR/0404040}


\bibitem{Rebbechi01}
\textbf{D\,Y Rebbechi}, \emph{Algorithmic properties of relatively hyperbolic
  groups}, PhD thesis, Rutgers Univ. (2001) \arxiv{math.GR/0302245}

\bibitem{SageevWiseTitsAlt}
\textbf{M Sageev}, \textbf{D\,T Wise}, \emph{The {T}its alternative for
  {$\CAT(0)$} cubical complexes}, \arxiv{math.GR/0405022}

\bibitem{Schroeder89}
\textbf{V Schroeder}, \emph{A cusp closing theorem}, Proc.\ Amer.\ Math.\ Soc.
  106 (1989) 797--802
  \MR{0957267}

\bibitem{Short91}
\textbf{H Short}, \emph{Quasiconvexity and a theorem of {H}owson's}, from:
  ``Group theory from a geometrical viewpoint \textup{(}{T}rieste,
  1990\textup{)}'', ({\'E} Ghys, A Haefliger, A Verjovsky, editors), World Sci.
  Publishing, River Edge, NJ (1991)  168--176
  \MR{1170365}

\bibitem{Wilson05}
\textbf{J\,M Wilson}, \emph{A {$\CAT(0)$} group with uncountably many distinct
  boundaries}, J.~Group Theory 8 (2005) 229--238 \MR{MR2126732}

\bibitem{WiseFigure8}
\textbf{D\,T Wise}, \emph{Subgroup separability of the figure $8$ knot group},
to appear in Topology

\bibitem{Wise96}
\textbf{D\,T Wise}, \emph{Non-positively curved squared complexes, aperiodic
  til\-ings, and non-residually finite groups}, PhD thesis, Princeton Univ.
  (1996)

\bibitem{XieTitsAltSquared04}
\textbf{X Xie}, \emph{Groups acting on {$\CAT(0)$} square complexes}, Geom.
  Dedicata 109 (2004) 59--88 \MR{MR2113187}

\bibitem{Xie4Manifold04}
\textbf{X Xie}, \emph{Tits alternative for closed real analytic $4$--manifolds
  of nonpositive curvature}, Topology Appl. 136 (2004) 87--121
  \MR{2023412}


\end{thebibliography}
\end{document}